\documentclass[reqno,a4paper,oneside]{amsart}
\usepackage{hyperref}
\usepackage{amssymb}
\usepackage{colortbl}
\usepackage{enumerate}
\usepackage{comment}
\allowdisplaybreaks

\begin{document}
	\title[Second-order conditions and Stability for the vche]
	{ Second-order optimality conditions and stability for optimal control problems governed by viscous Camassa-Holm equations}
	
	\author [C. T. Anh, N .H. H. Giang]
	{Cung The Anh and Nguyen Hai Ha Giang$^{\natural}$}

	\address{Cung The Anh \hfill\break
		Department of Mathematics, Hanoi National University of Education \hfill\break
		136 Xuan Thuy, Cau Giay, Hanoi, Vietnam}
	\email{anhctmath@hnue.edu.vn} 
	
	\address{Nguyen Hai Ha Giang \hfill\break
		Department of Mathematics, Hanoi National University of Education \hfill\break
		136 Xuan Thuy, Cau Giay, Hanoi, Vietnam}
	\email{giangnhh@hnue.edu.vn}

	\subjclass[2020]{49J20; 49K20; 76D55; 35Q35}
	\keywords{viscous Camassa-Holm equations; optimal control; second-order sufficient optimality condition;  stability; initial data.}
	
	\begin{abstract}  This work is a continuation of the previous one in [{\it Optimization} (2023)],  where the existence of optimal solutions and first-order necessary optimality conditions in both Pontryagin's maximum principle form and the variational form were proved for a distributed optimal control problem governed by the three-dimensional viscous Camassa-Holm equations in bounded domains with the cost functional of a quite general form and pointwise control constraints.  We will establish the second-order sufficient optimality conditions as well as the Lipschitz stability results of the control system with respect to perturbations of the initial data. 
	\end{abstract}
	\thanks{$^{\natural}$ Corresponding author: giangnhh@hnue.edu.vn}
	\maketitle
	
\numberwithin{equation}{section}
\newtheorem{theorem}{Theorem}[section]
\newtheorem{remark}{Remark}[section]
\newtheorem{definition}{Definition}[section]
\newtheorem{lemma}[theorem]{Lemma}
\newtheorem{corollary}[theorem]{Corollary}
\newtheorem{proposition}[theorem]{Proposition}
\newtheorem{example}[theorem]{Example}

\section{INTRODUCTION}
Let $\Omega$ be a smooth (at least $ C^3 $) bounded domain in $\mathbb{R}^{3}$ with boundary $\partial \Omega$, $ 0<T<+\infty $, and we denote the space-time cylinder by $Q=\Omega \times \left ( 0,T \right )$.  In this paper,  we consider the problem of finding a control $ h \in L^2(0,T; (L^2(\Omega))^3) $ and a state $ u \in W(0,T;H^3,V)  $ which minimize the cost function
\begin{equation} \label{subject function} \tag{\textbf{P}}
	  J(u,h) = \int_Q L(x,t,u(x,t),h(x,t))dxdt+\int_{\Omega}F(x,u(x, T))dx,
\end{equation}
where the control $ h $ and the state $ u $ satisfy the following 3D viscous Camassa-Holm equations (VCHE)
	\begin{equation}\label{state equation}
		\begin{cases}
		\partial_{t}\left ( u-\alpha ^{2} \Delta u \right )+\nu  \left (Au-\alpha ^{2} \Delta Au  \right ) + \nabla p\\
		\hspace{70pt}= u \times \left ( \nabla\times\left ( u-\alpha ^{2}\Delta u \right ) \right ) + h,  &(x,t) \in (\Omega \times (0,+\infty)),\\
		\nabla \cdot u = 0, &(x,t) \in (\Omega \times (0,+\infty)), \\
		u = Au = 0,   &(x,t) \in (\partial \Omega \times (0,+\infty)),  \\
		u(x,0) = u_0(x),  &\ x\in\Omega.
	\end{cases}
\end{equation} 
Here,  the parameter $ \alpha > 0$,  $\nu>0$ is the viscous constant,  $A$ is the Stokes operator,
$F: \Omega \times \mathbb{R}^3 \rightarrow \mathbb{R}$  and $L: Q\times\mathbb{R}^3\times\mathbb{R}^3 \rightarrow \mathbb{R}$ are given functions satisfying some certain conditions specified in Sections 2 and 3 below,  $ u $ is the fluid velocity  and the control $ h $ belongs to the following admissible set $\mathcal{H}$ of pointwise constraints
\begin{equation} \label{pointwise_constraint}
		h\in \mathcal{H}= \left \{ h\in (L^2(Q))^3 \ |\ h(x,t) \in \mathcal{M} \mbox{ for a.e. }(x,t)\in Q  \right \},
\end{equation}
where $\mathcal{M}$ is a non-empty closed convex set in $\mathbb{R}^3$.

A standard example for the choice of the cost function $ J $ is the quadratic function
\begin{equation*}
\begin{aligned}
	J(u,h) =\frac{\alpha_Q}{2} & \int_Q  |u(x,t)-u_d(x,t)|^2 dxdt +\frac{\alpha_T}{2} \int_{\Omega} |u(x,T)-u_T(x)|^2 \ dx \\&+ \frac{\gamma}{2} \int_Q|h(x,t)|^2 dxdt,  
	\end{aligned}
	\end{equation*}
where $u_d\in L^2(0,T;(L^2(\Omega))^3)$ and $u_T\in L^2(\Omega)$ denote the desired states of the system, and the term  $\frac{\gamma}{2} \int_Q |h(x,t)|^2 \ dxdt\; (\gamma>0)$ can be considered as the cost term.

It has been a longstanding problem in fluid dynamics to derive a model for the large scale motion of a fluid that averages or course-grains the small, computationally unresolvable, scales of the Navier-Stokes equations. The viscous Camassa-Holm equations (VCHE for short), also called Lagrangian averaged Navier-Stokes (LANS-$\alpha$) equations or Navier-Stokes-$\alpha$ equations, provide one such averaged model. The inviscid ($\nu=0$) case, known as the Lagrangian averaged Euler (LAE-$\alpha$) or Euler-$\alpha$ equations, was  introduced in \cite{Holm} as a natural mathematical generalization of the integrable inviscid 1D Camassa-Holm equation discovered in \cite{Camassa} through a variational formulation. In \cite{Foias2002}, the authors first added viscous dissipation to the equations, they argued on physical grounds that the momentum $u+\alpha^2Au$ rather than the velocity $u$, need to be diffused.  The main reason of studying the VCHE comes from the need of approximating problems relating to turbulent flows, because this kind of models preserves properties of transport for circulation and vorticity dynamics of the Navier-Stokes equations.  For a complete description of the physical significance of the VCHE, namely in turbulence theory, and their developments, we refer \cite{Foias2002,  Marsden2001} and references therein.

In the past years, the existence and long-time behavior of solutions to the VCHE have attracted the attention of many mathematicians. In bounded domains with Dirichlet or periodic boundary conditions, there are many results on the existence of solutions and existence of global attractors for VCHE, see e.g. \cite{Coutand, Foias2002, Ilyn2003, Kim2006, Marsden2001, Vishik2007} and references therein. The time decay rates of solutions on the whole space were extensively investigated in \cite{AT2017, Bjorland2008,  BjorlandSchonbek, Gao2022, Zhao}. We also refer the interested reader to \cite{Albanez2016, AnhBach, Korn2009} for recent results on the data assimilation to the VCHE and to \cite{Gao2016, Mitra} for the controllability of the  VCHE. 

On the other hand,  the optimal control plays an important role in modern control theories, and has a wider application in modern engineering. Modern optimal control theories and applied models are both represented by ODEs, which have been developed perfectly. With the development and application of technology, it is necessary to solve the problems of optimal control theory for PDEs. Especially, considerable progress has been made in mathematical analyses and computations of optimal control problems for viscous flows. The optimal control problems for the Navier-Stokes kind equations have been studied extensively during the past years,  see e.g. \cite{Abergel, Fattorini, Gunzburger, Gunzburger1999, Hinze2002, HK2001, Sritharan,  Tr2005, Wachsmuth1, Wachsmuth2006, Wang2002} and references therein.  

Let us briefly review some existing results on optimal control of viscous Camassa-Holm equations (VCHE).  In the paper \cite{Tian2009},  the authors proved the existence of optimal solutions for VCHE in the one-dimensional case with a distributed control, but no optimality conditions were given there.  Later,  in the paper \cite{AS2019} the authors proved the existence of optimal solutions and established the first-order necessary as well as the second-order sufficient optimality conditions for an optimal control problem of the three-dimensional VCHE in bounded domains with a standard quadratic objective functional.  See also the related work \cite{MOV2021},  where the solvability of the optimal control problem,  the first-order optimality conditions and convergence,  as the parameter $\alpha \to 0^+$,  of the optimality system of the optimal control problem associated to the VCHE  were investigated.  The time optimal control of VCHE was studied in \cite{SonThuy}.  In the  recent work \cite{AGN2023},  we proved the existence of optimal solutions and derived the first-order necessary optimality conditions in both Pontryagin's maximum principle form and the variational form for a distributed optimal control problem governed by the three-dimensional VCHE  in bounded domains with a quite general class of cost functionals and  pointwise control constraints. These results can be seen in some sense as the natural generalizations of the corresponding ones in previous work \cite{AS2019}.   However,  as far as we know,  the second-order optimality conditions and the stability results  to the optimal control problem for the VCHE, especially in the case of pointwise control constraints,  have not been studied before.  This is the main motivation of our present paper.

We address two main issues in this paper.  First,  we prove the second-order sufficient optimality conditions of problem \textbf{(P)} in both forms with a quite general objective functional, and in this sense the second-order optimality conditions obtained here and the first-order ones in the previous paper \cite{AGN2023} are natural generalizations of those  in \cite{AS2019}.  Second, we investigate the stability of selected local solutions of the optimal control problem \textbf{(P)} with respect to a perturbation of the initial data. To the best of our knowledge,  this is the first result on the stability of the optimal solutions of the VCHE with respect to perturbations.  It is noticed that in most of papers dealing with the stability,  perturbations appeared in the differential equation, in its boundary condition, in the objective functional, or in inequality constraints (see e.g.  \cite{CT2016,  Corella2023, Khanh2024} for recent works on semilinear parabolic optimal control problems). As mentioned in \cite{CT2022},  handling  perturbations of the initial data is more complicated for a nonlinear state equation, in particular since bounded initial data are needed to have a differentiable control-to-state mapping.  Beside the similarity of our approach and techniques to those in \cite{AS2019,CT2022}, due to the generalized form of the cost functional and the complexity of nonlinear terms in the VCHE, we have met some new difficulties in the proof of main results. To overcome this,  we particularly need to show the twice differentiability of the control-to-state mapping and exploit the contradiction arguments to derive the second-order optimality conditions,  and using the sufficient conditions to show the stability results. The methods used in this paper can be applied to optimal control problems of other PDEs with a similar structure.

The paper is organized as follows. In Section 2, we derive some auxiliary results concerning both the states equations and the objective functional, which are frequently used later. After reformulating necessary assumptions on the functions in the objective functional, we establish the second-order sufficient optimality conditions in Section 3. Under the mentioned optimality conditions imposed on local optimal solutions, the associated stability analysis is performed in Section 4. The main result of this section is Theorem \ref{3rd_theorem_stability} on the Lipschitz stability of local solutions with respect to a perturbation in the initial data.  
\section{PRELIMINARIES AND AUXILIARY RESULTS}
\subsection{Function spaces and inequalities for the nonlinear terms} 
Throughout this paper, we shall denote by $ (\cdot,\cdot) $ and $ |\cdot|  $, the scalar product and the associated norm in $ (L^2(\Omega))^3 $, respectively; and by $ (\nabla u,\nabla v) $ the scalar product in $ (L^2(\Omega))^3 $ of the gradients of $ u $ and $ v $. Also, we define $ ((u,v)) = (\nabla u,\nabla v) $ for $ u,v \in (H_{0}^{1}(\Omega))^{3} $ is the scalar product in $(H_{0}^{1}(\Omega))^{3} $; and its associated norm $ \left\| \cdot  \right\| $.

Hereafter we shall use the following function spaces
\begin{align*}
	&\mathbb{L}^2(\Omega) := (L^2(\Omega))^3; \\
	&\mathbb{L}^2(Q) := (L^2(Q))^3; \\
	&\mathbb{H}^k(\Omega) := (H^k(\Omega))^3;\\
	& \mathcal{V} = \{ u \in (C_0^{\infty}(\Omega))^3: \nabla \cdot u = 0 \}; \\
	& H := \text{the closure of }  \mathcal{V} \ \text{in} \ \mathbb{L}^2(\Omega)= \{ u \in \mathbb{L}^2(\Omega)\ |\ \nabla \cdot u = 0 \ \text{and} \ u \cdot \vec{n} = 0 \ \text{on} \ \partial \Omega \};\\
	& V := \text{the closure of }  \mathcal{V} \ \text{in} \ \mathbb{H}_0^1(\Omega) = \{ u \in \mathbb{H}_0^1(\Omega) \ |\ \nabla \cdot u = 0 \};\\
	& \mathbb{H}_{\sigma}^k := \text{the closure of }  \mathcal{V} \ \text{in} \ \mathbb{H}_{0}^k(\Omega); \\
	& D(A) = (H^2(\Omega))^3 \cap V;\ V' \ \text{and} \ D(A)' \ \text{are the dual space of} \ V \text{and} \ D(A).
\end{align*}
Moreover, we introduce the space of functions $ u $ whose time derivatives $\frac{du}{dt}$ exist as abstract functions
\[ W(0,T;H^3,V) = \{ u \in L^2(0,T;\mathbb{H}_{\sigma}^3(\Omega)) \ |\ \frac{du}{dt} \in L^2(0,T;V) \}. \]
It is well-known that the embeddings
\[ D(A) \hookrightarrow V \hookrightarrow H \]
are compact and each space is dense in another space. Let us denote by $ \mathcal{P}: \mathbb{L}^2(\Omega) \rightarrow H $ the Leray projection on $ H $, and by $ A: D(A) \rightarrow H, A = -\mathcal{P}\Delta $ the Stokes operator with domain $ D(A) $. Since $ \partial \Omega $ is smooth, $ |Au| $ defines in $ D(A) $ a norm which is equivalent to the $ H^2(\Omega) $-norm (see   \cite[Proposition 4.7]{Constantin}), and therefore $ D(A) $ is a Hilbert space with the scalar product $ (u,v)_{D(A)} = (Au, Av) $.

Similarly, when $ Au = 0 $ on $ \partial \Omega $, the operator $ A^2 $ can be extended continuously to be defined on $ D(A) $ with value in $ D(A)'  $ as following
\[ \left \langle A^2u,v \right \rangle_{V', V} = (Au,Av) \ \forall u, v \in D(A). \]
Following the notations for the Navier-Stokes equations we denote the trilinear form 
\[ b(u,v,w) = \sum_{i,j = 1}^{3} \int_{\Omega} u_i \frac{\partial v_j}{\partial x_i} w_j dx, \]
whenever the integral makes sense.
According to \cite{AS2019}, we have the following properties for $ b(u,v,w) $
\begin{align*}
	b(u,v,w) &= -b(u,w,v), \ \forall(u, v, w)  \in V \times V \times V, \nonumber \\
	|\left \langle u \cdot \nabla v,w \right \rangle | &= | b(u,v,w) | \leq 
	C \begin{cases} 
		\left\| u \right\| |v| |Aw|, \ \forall (u, v, w)  \in V \times H \times D(A),\\ 
		|u| \left\| v \right\| |Aw|, \ \forall  (u, v, w)  \in H \times V \times D(A),\\
	\end{cases}\\  
	|\left \langle v\cdot \nabla u^T,w \right \rangle | &= |b(w,u,v)| \leq C \left\| u \right\| |v| |Aw|, \  \forall  (u, v, w)  \in V \times H \times D(A).
\end{align*}
Furthermore, we can define a continuous bilinear operator $ \widetilde{B}$ from $ V \times V $ into $ V'$ by \[ \left \langle \widetilde{B}(u,v),w \right \rangle _{V',V} = \widetilde{b}(u,v,w), \]
where $ \widetilde{b}(u,v,w) = b(u,v,w) - b(w,v,u). $
As stated in \cite{Foias2002},  the following results hold for $ \widetilde{B}$
\begin{align*}
	(i)\ & \left \langle \widetilde{B}(u,v),w \right \rangle_{V',V} = - \left \langle \widetilde{B}(w,v),u \right \rangle_{V',V}, \\
	(ii)\ &  \left \langle \widetilde{B}(u,v),u \right \rangle_{V',V} \equiv 0, \\
	(ii)\ &  \left \langle \widetilde{B}(u,v),w \right \rangle_{V',V}  \leq C \vert u \vert^{1/2} \Vert u \Vert^{1/2} \Vert v \Vert \Vert w \Vert,\\
	(iii)\ & \left \langle \widetilde{B}(u,v),w \right \rangle_{V',V}  \leq C \Vert u \Vert \Vert v \Vert \vert w \vert^{1/2} \Vert w \Vert^{1/2}, 
\end{align*}
for every  $ u,v,w \in V $, and  $ C $ is a positive constant depending on only $\Omega $ .
\subsection{Existence and uniqueness of weak solution to the viscous Camassa-Holm equations (VCHE)}
\begin{definition}\label{VCHE_weaksolution}
	Let $ h \in L^2(0,T;\mathbb{L}^2(\Omega)) $ and let $ T > 0 $. A function \[ u \in L^\infty(0,T;\mathbb{H}_{\sigma}^2(\Omega)) \cap L^2(0,T;\mathbb{H}_{\sigma}^3(\Omega)) \] with $\displaystyle \frac{du}{dt} \in L^2(0,T;V)  $ is said to be a weak solution to problem \eqref{state equation} on the interval $ (0,T)  $ if it satisfies
	\begin{align} \label{weak solution}
		 \left \langle \partial_{t}u + \alpha^2 A\partial_{t}u, w \right \rangle_{V', V} &+ \nu \left \langle A(u + \alpha^2 Au), w \right \rangle_{V', V} \nonumber \\
		& + \left \langle \widetilde{B}(u, u+\alpha^2 Au), w \right \rangle_{V', V} = \langle h,w \rangle_{V',V}
	\end{align}
	for every $ w \in V $ and for almost every $ t \in [0,T] $. Moreover, $ u(0) = u_0 $ in $ D(A) $.
\end{definition}
Here, equation \eqref{weak solution} is understood in the following sense: for a.e. $ t_0, t \in [0,T] $ and for all $ w \in V $, we have
\begin{align*}
	(u(t) &+ \alpha^2 Au(t),w) - (u(t_0) + \alpha^2 Au(t_0),w) + \nu \int_{t_0}^{t} (u(s) + \alpha^2 Au(s),Aw) ds \\
	&+ \int_{t_0}^{t} \left \langle \widetilde{B}(u(s),u(s) + \alpha^2 Au(s)),w \right \rangle {_{V', V}} ds = \int_{t_0}^{t} (h,w) ds.
\end{align*}
The following theorem is proved by using the arguments in \cite{BjorlandSchonbek} (or \cite{Coutand}).
\begin{theorem}
	For $ h \in L^2(0,T;\mathbb{L}^2(\Omega)), T >0 $ and $ u_0 \in D(A) $ given, there exists a unique weak solution to problem \eqref{state equation} on the interval $ (0, T) $ in the sense of Definition \ref{VCHE_weaksolution}. Moreover, there exists a constant $ C $ such that the function $ u $ satisfies the following estimate for all $ t \in [0, T] $,
	\begin{align*}
		|u(t)|^2 &+\alpha^2 \|u(t)\|^2+\nu \int_{0}^{t}(\|u(s)\|^2+\alpha^2 |Au(s)|^2)ds \\ & \leq |u_0|^2 + \alpha^2 \|u_0\|^2+ C \|h\|^2_{L^2(0,T;\mathbb{L}^2(\Omega))}. 
	\end{align*}
\end{theorem}
\begin{remark}
	{\rm The boundedness of $ u $ in $ (L^{\infty}(Q))^3 $ can be shown if $ h $ is bounded in $ L^2(0,T;\mathbb{L}^2(\Omega)) $. }
\end{remark}

\subsection{The control-to-state mapping}
Let us study the behavior of the mapping: right-hand side $ \mapsto $ solution, the so-called \textit{control-to-state mapping} and its Fréchet differentiability. This subsection follows the general lines of the approach used in \cite{Wachsmuth2006}.
\begin{definition} \label{Smapping}
	Consider the system \eqref{state equation}. The mapping $ h \mapsto u $, where $ u $ is the weak solution of \eqref{state equation} with fixed initial value $ u_0$, is denoted by $ S $, i.e. $ u = S(h) $.
\end{definition}

\begin{lemma}\label{LipchitzS}
 Let $u_1,u_2 \in W(0,T;H^3,V)$ be the weak solutions of \eqref{state equation} with the right-hand sides of the first equations equal to $h_1$ and $h_2$ in $\mathbb{L}^2(Q)$ and the initial datas equal to $u_{01}$ and $u_{02}$, respectively. Then there exists a constant $C$ such that the following estimate holds
\begin{align}
	\| u_1 - u_2 \|^2_{W(0,T;H^3,V)} \leq C\left(\|u_{01} - u_{02}\|^2 + \| h_1 - h_2\|_{\mathbb{L}^2(Q)}^2 \right).
\end{align}
\end{lemma}
\begin{proof}
	Putting $ \delta u = u_1 - u_2 $, thus $ \delta u $ is the solution of the following equations
	\begin{equation} \label{estimate_deltau}
		\begin{cases}
			\partial_t(\delta u +\alpha^2 A\delta u)+\nu A(\delta u+\alpha^2 A\delta u)
			+\widetilde{B}(u_1,u_1+\alpha^2 A u_1)\\ 
			\hspace{110pt}- \widetilde{B}(u_2,u_2+\alpha^2 Au_2) = h_1 - h_2 \ \text{in} \ L^2(0,T;V'),\\
			\delta u(0)=u_{0_1} - u_{0_2},
		\end{cases}
	\end{equation}
	where  \[ \widetilde{B}(u_1,u_1+\alpha^2 A u_1) - \widetilde{B}(u_2,u_2+\alpha^2 Au_2) = \widetilde{B}(\delta u, u_1+\alpha^2 A u_1) + \widetilde{B}(u_2, \delta u+ \alpha^2 A \delta u) . \]
	Taking the inner product of \eqref{estimate_deltau} with $ \delta u $ yields the following equality
\[  \frac{1}{2} \frac{d}{dt}(|\delta u|^2+\alpha^2 \left \| \delta u \right \|^2) +\nu (\left \| \delta u \right \|^2 + \alpha^2 |A\delta u|^2)+\widetilde{b}(u_2, \delta u+\alpha^2 A\delta u,\delta u) = (h_1 - h_2,\delta u) . \]
	Applying Young's inequality we obtain
	\[  \vert \widetilde{b}(u_2, +\alpha^2 A\delta u,\delta u) \vert \leq C \left \|\delta u \right \| ^2 +\frac{\nu \alpha^2}{2} |A\delta u|^2, \]
	and
	\[ \vert (h_1 - h_2,\delta u) \vert \leq C_1 \vert h_1 - h_2 \vert^2 +\frac{\nu}{2} \left \|  \delta u \right \| ^2.  \]
	It follows from these estimates that
	\begin{align*}
		\frac{1}{2}\frac{d}{dt}(|\delta u|^2&+\alpha^2 \left \| \delta u\right \|^2) +\frac{\nu}{2} (\left \| \delta u \right \|^2 + \alpha^2 |A\delta u|^2) \leq C \Vert \delta u \Vert ^2 + C_1 \vert h_1 - h_2 \vert^2\\
		&\leq \frac{C}{\alpha^2}(|\delta u|^2 + \alpha^2 \Vert \delta u \Vert^2) + C_1 \vert h_1 - h_2 \vert^2.
	\end{align*}
	Thanks to Gronwall's inequality we get
	\begin{align*}
		|\delta u(t)|^2 &+\alpha^2 \left \| \delta u (t) \right \|^2 + \frac{\nu}{2} \int_{0}^{t} \left(\left\| \delta u(s) \right\|^2 +\alpha^2 |A\delta u(s)|^2 \right) ds 
		\\ & \leq \left ( \alpha^2 \Vert \delta u(0)\Vert ^2 + |\delta u(0)|^2 \right ) e^{Ct} + \int_{0}^{t} e^{C(t-s)} \vert (h_1 - h_2)(s) \vert^2ds. 
	\end{align*}
	Thus, $ \Vert \delta u \Vert^2_{L^2(0,T;D(A))} \leq C\left(\|u_{01} - u_{02}\|^2 + \| h_1 - h_2\|_{\mathbb{L}^2(Q)}^2 \right)$. 
	
	We now estimate $ \Vert \delta u \Vert_{L^2(0,T;\mathbb{H}_{\sigma}^3(\Omega))} $. Since $ A \delta u$  is divergence-free, we can take the inner product of \eqref{estimate_deltau} with $ A \delta u $ to obtain
	\begin{align}
		\frac{1}{2} \frac{d}{dt} \left(\Vert \delta u \Vert^2 + \alpha^2 \vert A\delta u \vert^2 \right) &+ \nu \left( \vert A\delta u \vert ^2 + \alpha^2 \vert \nabla A \delta u \vert^2 \right) + \widetilde{b} (\delta u, u_1 + \alpha^2 A u_1, A \delta u) \nonumber \\
		&+\widetilde{b} (u_2, \delta u + \alpha^2 A \delta u, A \delta u) = (h_1-h_2, A\delta u). \label{innerproductAdeltau}
	\end{align}
Then, since $ u_1, u_2 \in W(0,T;H^3,V) $ and $ \delta u \in L^2(0,T;D(A)) $, we approach the nonlinear terms of \eqref{innerproductAdeltau} in the following way
\begin{align*}
	\vert (h_1 - h_2, A\delta u) \vert &\leq \vert h_1 - h_2 \vert \vert A\delta u \vert \leq \vert h_1 - h_2 \vert^2 + \frac{\nu}{4} \vert A \delta u \vert^2, \\
	\vert \widetilde{b}(\delta u, u_1 + \alpha^2 A u_1, A \delta u) \vert &\leq \vert \widetilde{b}(\delta u, u_1, A\delta u) \vert + \alpha^2 \vert \widetilde{b} (\delta u, Au_1, A\delta u) \vert \\
	&\leq C\Vert \delta u \Vert_{L^4} \vert \nabla u_1 \vert \Vert A\delta u \Vert _{L^4} + C \Vert \delta u \Vert_{L^4} \vert \nabla Au_1 \vert \Vert A\delta u \Vert_{L^4}\\
	&\leq C \vert A\delta u \vert ^{1/4} \vert \nabla A \delta u \vert ^{3/4} \Vert \delta u \Vert \\
	&\leq \frac{\nu \alpha^2 }{4} \vert \nabla A\delta u \vert^2 + C \vert A \delta u \vert^2 + C.
\end{align*}
To bound the term $ \widetilde{b}(u_2, \delta u + \alpha^2 A \delta u, A \delta u) $, start with the definition of the operator and Young's inequality
\begin{align*}
		\vert \widetilde{b}(u_2, \delta u + \alpha^2 A \delta u, A \delta u) \vert &\leq \vert \widetilde{b}(u_2, \delta u, A \delta u) \vert + \alpha^2 \vert \widetilde{b}(u_2, A\delta u, A\delta u) \vert \\
		&\leq \vert b(u_2, A\delta u, \delta u) \vert + \vert b (A\delta u, u_2, \delta u) \vert + \alpha^2 \vert b(A\delta u,u_2, \delta u ) \vert\\
		&\leq C \Vert u_2 \Vert_{L^4} \vert \nabla A\delta u \vert \Vert \delta u \Vert_{L^4} + C \Vert A\delta u \Vert_{L^4} \vert \nabla u_2 \vert \Vert \delta u \Vert_{L^4} \\
		&\hspace{20pt}+ C \Vert u_2 \Vert_{L^4} \vert \nabla A\delta u \vert \Vert A\delta u \Vert_{L^4}\\
		&\leq C \vert \nabla A\delta u \vert \Vert \delta u \Vert + C \vert A\delta u \vert^{1/4} \vert \nabla A\delta u \vert^{3/4} \Vert \delta u \Vert \\
		&\hspace{20pt}+ C \Vert u_2 \Vert_{L^4} \vert \nabla A\delta u \vert \Vert A\delta u \Vert_{L^4} \\
		&\leq \frac{\nu \alpha^2}{8} \vert \nabla A\delta u \vert^2 + C \vert A\delta u \vert^2 + C + C \vert A\delta u \vert^{1/4} \vert \nabla A\delta u \vert^{7/4}\\
		&\leq \frac{\nu \alpha^2}{8} \vert \nabla A\delta u \vert^2 + C \vert A\delta u \vert^2 + C \\
		&\hspace{20pt}+ C \vert A\delta u \vert^{1/12} \vert A\delta u  \vert^{1/6} \vert \nabla A\delta u \vert^{7/4}\\
		&\leq \frac{\nu \alpha^2}{8} \vert \nabla A\delta u \vert^2 + C \vert A\delta u \vert^2 + C + C \Vert \delta u \Vert^{1/12} \vert \nabla A\delta u \vert^{23/12} \\
		&\leq \frac{\nu \alpha^2}{4} \vert \nabla A\delta u \vert^2 + C \vert A\delta u \vert^2 + C.
\end{align*}
Putting all these bounds together yields
\begin{align}
	\vert \widetilde{b}(\delta u&, u_1+\alpha^2 Au_1, A\delta u)\vert + \vert \widetilde{b}(u_2, \delta u+\alpha^2 A\delta u, A\delta u) \vert + \vert (h_1 - h_2, A\delta u) \vert \nonumber \\
	&\leq \frac{\nu \alpha^2}{2} \vert \nabla A\delta u \vert^2 + C \vert A\delta u \vert^2 + \frac{1}{\nu} \vert h_1 - h_2 \vert^2 + C \Vert \delta u \Vert^2. \label{bAdeltau}
\end{align}
Using \eqref{bAdeltau} we deduce that
\begin{align}
	\frac{1}{2} \frac{d}{dt} & \left(\Vert \delta u \Vert^2 + \alpha^2 \vert A \delta u \vert^2 \right) + \nu \left(\vert A\delta u \vert^2 + \alpha^2 \vert \nabla A \delta u \vert^2 \right) \nonumber \\
	& \leq \frac{\nu \alpha^2 }{2} \vert \nabla A \delta u \vert^2 + C \vert A \delta u \vert^2 + C \Vert \delta u \Vert^2 + \frac{1}{\nu} \vert h_1 - h_2 \vert^2. \label{estimateAdeltau}
\end{align}
For $ t>0 $, integrate \eqref{estimateAdeltau} from $0 $ to $ t $ to see
\begin{align*}
	\Vert \delta u(t) \Vert ^2 &+ \alpha^2 \vert A \delta u(t) \vert^2 + \frac{\nu}{2} \int_{0}^{t} \vert A \delta u(s) \vert^2 ds + \frac{\nu \alpha^2}{2} \int_{0}^{t} \vert \nabla A \delta u (s) \vert^2 ds \\
	&\leq \vert A \delta u (0) \vert^2 + \Vert \delta u (0) \Vert^2 + \frac{1}{\nu} \int_{0}^{t} \vert (h_1 - h_2)(s) \vert^2 ds.
\end{align*}
Combining with $ \delta u(0) \in D(A) $, we conclude that
\[ \Vert \delta u \Vert^2_{L^2(0,T;\mathbb{H}_{\sigma}^3(\Omega))} \leq C\left(\Vert u_{0_1} - u_{0_2} \Vert^2 + \| h_1 - h_2\|_{\mathbb{L}^2(Q)}^2 \right). \]
The remaining norm $ \Vert \delta u_t \Vert_{L^2(0,T;V)} $ is handled in a similar way. Multiplying \eqref{estimate_deltau} with $ \delta u_t \in L^2(0,T;V) $ shows
\begin{align*}
	\vert \delta u_t \vert^2 &+ \alpha^2 \Vert \delta u_t \Vert^2 + \nu \int_{0}^{T} ((\delta u(s), \delta u_t(s) )) ds + \nu \alpha^2 \int_{0}^{T} \left(A\delta u(s), A\delta u_t(s)\right) ds\\
	& + \widetilde{b}(\delta u,u_1+\alpha^2Au_1, \delta u_t) + \widetilde{b}(u_2, \delta u + \alpha^2 A\delta u, \delta u_t) = \left \langle h_1 - h_2, \delta u_t \right \rangle_{V', V}.
\end{align*}
Again, owing to the fact that $ u_1, u_2 \in W(0,T;H^3,V) $, and Young's inequality,  we have
\begin{align*}
	\left \vert \left \langle h_1 - h_2, \delta u_t \right \rangle_{V', V} \right \vert &= \int_{0}^{T} \left(h_1(s) - h_2(s), \delta u_t(s)\right)ds \leq C  \vert h_1 - h_2 \vert^2 + \frac{1}{2} \vert \delta u_t \vert^2,\\
	\int_{0}^{T} ((\delta u(s), \delta u_t(s))) ds &= \Vert \delta u(T) \Vert^2 - \Vert \delta u(0) \Vert^2,\\
	\int_{0}^{T} \left(A\delta u(s), A\delta u_t(s)\right) ds &= \vert A\delta u(T) \vert^2 - \vert\delta u(0) \vert^2,
\end{align*}
and
\begin{equation*}
	\vert \widetilde{b}(\delta u, u_1 + \alpha^2 Au_1, \delta u_t) \vert + \vert \widetilde{b}(u_2, \delta u + \alpha^2 A\delta u, \delta u_t) \vert
	\leq C \Vert \delta u_t \Vert^2 + \frac{1}{2} \vert\delta u_t \vert^2 + C.
\end{equation*}
Putting these inequalities together,
\[ C\Vert\delta u\Vert^2 + \nu \left(\Vert\delta u(T)\Vert^2 - \Vert \delta u(0)\Vert^2\right) + \nu \alpha^2 \left(\vert A\delta u(T) \vert^2 - \vert A\delta u(0) \vert^2 \right) \leq C\vert h_1 - h_2 \vert^2 + C. \]
This relation gives an estimate on $ \delta u_t $
\[ C\Vert \delta u_t \Vert^2 + \nu \Vert \delta u(T) \Vert^2 + \nu \alpha^2 \vert A\delta u(T) \vert^2 \leq C\vert h_1 - h_2 \vert^2 + \nu \Vert \delta u(0) \Vert^2 + \nu \alpha^2 \vert A\delta u(0) \vert^2, \]
or 
\[ \Vert \delta u_t \Vert^2_{L^2(0,T;V)} \leq C \left(\Vert u_{0_1} - u_{0_2} \Vert^2 + \| h_1 - h_2\|_{\mathbb{L}^2(Q)}^2 \right). \]
The proof is complete.
\end{proof}
\begin{theorem}\label{1stS}
	The control-to-state mapping is Fréchet differentiable as mapping from $ L^2(0,T;\mathbb{L}^2(\Omega)) $ to $ W(0,T;H^3,V) $. The derivative at $ \hat{h} \in L^2(0,T;\mathbb{L}^2(\Omega)) $ in direction $ k \in L^2(0,T;\mathbb{L}^2(\Omega)) $ is given by $ S'(\hat{h})h = u $, where $ u $ is the weak solution of 
		\begin{equation} \label{linearized eq}
		\begin{cases}
			\partial_{t} (u + \alpha^{2} Au) + \nu A(u + \alpha^{2} Au) + \widetilde{B}'(\hat{u}, \hat{u}+\alpha^2 A \hat{u})u = k \ \text{in} \ L^2(0,T;V'),\\
			u(0)=0 \ \text{in} \ D(A),
		\end{cases}
	\end{equation}
with $ \hat{u} = S(\hat{h}) $ and $ \widetilde{B}'(\hat{u}, \hat{u}+\alpha^2 A \hat{u})u = \widetilde{B}(u, \hat{u} + \alpha^{2} A \hat{u}) + \widetilde{B}(\hat{u}, u + \alpha^{2}A u) $.
\end{theorem}
\begin{proof}
	Define $ u = S(\hat{h}+k) $. We find that the difference $ d := u - \hat{u} $ is the weak solution of
	\begin{equation} \label{u-hatu1}
		\begin{cases}
			\partial_{t} (d + \alpha^{2} Ad) + \nu A(d + \alpha^{2} Ad) + \widetilde{B}(u, u +\alpha^2 Au) - \widetilde{B}(\hat{u}, \hat{u}+\alpha^2 A \hat{u}) = k ,\\
			d(0)=0 \ \text{in} \ D(A).
		\end{cases}
	\end{equation}
However, since \[ \widetilde{B}(u, u +\alpha^2 Au) - \widetilde{B}(\hat{u}, \hat{u}+\alpha^2 A \hat{u})  =  \widetilde{B}'(\hat{u}, \hat{u} + \alpha^2 A \hat{u})d + \widetilde{B}(d, d+\alpha^2 Ad), \] equations \eqref{u-hatu1} can be transformed into
	\begin{equation*}
	\begin{cases}
		\partial_{t} (d + \alpha^{2} Ad) + \nu A(d + \alpha^{2} Ad) + \widetilde{B}'(\hat{u}, \hat{u} + \alpha^2 A \hat{u})d = k -\widetilde{B}(d, d+\alpha^2 Ad)  ,\\
		d(0)=0 \ \text{in} \ D(A).
	\end{cases}
	\end{equation*}
We split $ d $ into $ d = z+r $, where $ z $ and $ r $, respectively, are the weak solutions of these systems
	\begin{equation*}
	\begin{cases}
		\partial_{t} (z + \alpha^{2} Az) + \nu A(z + \alpha^{2} Az) + \widetilde{B}'(\hat{u}, \hat{u} + \alpha^2 A \hat{u})z = k  ,\\
		z(0)=0,
	\end{cases}
\end{equation*}
and
	\begin{equation*}
	\begin{cases}
		\partial_{t} (r + \alpha^{2} Ar) + \nu A(r + \alpha^{2} Ar) + \widetilde{B}'(\hat{u}, \hat{u} + \alpha^2 A \hat{u})r = -\widetilde{B}(d, d+\alpha^2 Ad)  ,\\
		r(0)=0.
	\end{cases}
\end{equation*}
By similar estimates as those in Lemma \ref{LipchitzS}, we obtain the boundedness of $ \Vert z \Vert_{W(0,T;H^3,V)} $ and $ \Vert r \Vert_{W(0,T;H^3,V)} $,
\[ \Vert z \Vert_{W(0,T;H^3,V)} \leq C \Vert k \Vert_{L^2(0,T;V')}, \]
\[ \Vert r \Vert_{W(0,T;H^3,V)} \leq C \Vert \widetilde{B}(d,d+\alpha^2Ad)\Vert_{L^2(0,T;V')}. \]
Meanwhile, we can derive
\begin{align*}
	\left | \left \langle \widetilde{B}(d, d+\alpha^2 Ad), v \right \rangle_{V',V} \right | &= |\widetilde{b}(v, d+\alpha^2Ad, d)| \\
	& \leq C\Vert v \Vert \Vert d \Vert^{1/2} |d+\alpha^2 Ad|^{3/2} \\
	& \leq C\Vert d\Vert^2 + \nu \alpha^2 |Ad|^2\\
	& \leq C\Vert v\Vert \Vert d \Vert^2.
\end{align*}
Hence, it holds $ \Vert r \Vert_{W(0,T;H^3,V)} \leq C \Vert d \Vert^2 = \Vert S(\hat{h}+k) - S(\hat{h}) \Vert^2_{W(0,T;H^3,V)} $.\\
Taking into account the property of the mapping $ S $ given by Lemma \ref{LipchitzS},
the inequality $ \Vert r \Vert_{W(0,T;H^3,V)} \leq C  \Vert k \Vert^2_{L^2(0,T;\mathbb{L}^2(\Omega))}$ follows. Thus, we have shown that 
\[ \frac{\Vert u - \hat{u} - z \Vert_{W(0,T;H^3,V)}}{\Vert k \Vert_{L^2(0,T;\mathbb{L}^2(\Omega))}} \rightarrow 0 \ \text{as} \ \Vert k \Vert_{L^2(0,T;\mathbb{L}^2(\Omega))} \rightarrow 0,\]
which implies that the function $ z $, denoted by $ z = S'(\hat{h})k $, will be the Fréchet derivative of $ S' $ at $ \hat{h} $ in direction $ k $ and $ \Vert S'(\hat{h})k \Vert_{W(0,T;H^3,V)} \leq C \| k\|_{\mathbb{L}^2(Q)}^2$. This completes the proof.
\end{proof}

\begin{theorem}
	The control-to-state mapping is twice continuously differentiable as mapping from $ L^2(0,T;\mathbb{L}^2(\Omega)) $ to $ W(0,T;H^3,V) $. The derivative at $ \hat{h} \in L^2(0,T;\mathbb{L}^2(\Omega))$ in directions $ h_1, h_2 \in L^2(0,T;\mathbb{L}^2(\Omega))$ is given by $ S''(\hat{h})[h_1,h_2] = u $, where $ u $ is a weak solution of
	 	\begin{equation} \label{2nd_S}
	 	\begin{cases}
	 		\partial_{t} (u + \alpha^{2} Au) + \nu A(u + \alpha^{2} Au) + \widetilde{B}'(\hat{u}, \hat{u}+\alpha^2 A \hat{u})u = -\widetilde{B}''(\hat{u},\hat{u}+\alpha^2 A\hat{u})(u_1,u_2) ,\\
	 		u(0)=0 \ \text{in} \ D(A),
	 	\end{cases}
	 \end{equation}
 with $ \hat{u} = S(\hat{h}) $ and $ u_i = S'(\hat{h})h_i, i = 1,2 $.\\
 Here, the term $ \widetilde{B}''(\hat{u},\hat{u}+\alpha^2 A\hat{u})(u_1,u_2)  $ is decomposed as
 \[ \widetilde{B}''(\hat{u},\hat{u}+\alpha^2 A\hat{u})(u_1,u_2) = \widetilde{B}(u_1, u_2 + \alpha^2 Au_2) + \widetilde{B}(u_2, u_1 + \alpha^2 A u_1). \]
\end{theorem}
\begin{proof}
Let be given $ h_1, h_2 \in L^2(0,T;\mathbb{L}^2(\Omega)). $\\
Putting $ u_1 = S'(\hat{h})h_1 $ and $ \tilde{u} = S'(\hat{h}+h_2)h_1 $. By Theorem \ref{1stS}, $\tilde{u} $ is the weak solution of 
	\begin{equation*}
	\begin{cases}
		\partial_{t} (\tilde{u} + \alpha^{2} A\tilde{u}) + \nu A(\tilde{u} + \alpha^{2} A\tilde{u}) \\ \hspace{20.5pt} + \widetilde{B}' \left[ S(\hat{h}+h_2), S(\hat{h}+h_2)+\alpha^2 A \, S(\hat{h}+h_2) \right] \tilde{u} = h_1 \ \text{in} \ L^2(0,T;V'),\\
		\tilde{u}(0)=0 \ \text{in} \ D(A).
	\end{cases}
\end{equation*}
We want to write $ \tilde{u} $ as $ S'(\hat{h})f $ with some function $ f $. Using the Fréchet differentiability of $ S $, we find $ S(\hat{h}+h_2) = S(\hat{h}) + S'(\hat{h})h_2+r_1(h_2) $ with some remainder term $ r_1 $ satisfying 
\[ \frac{\Vert r_1(h_2) \Vert_{W(0,T;H^3,V)}}{\Vert h_2 \Vert_{L^2(0,T;\mathbb{L}^2(\Omega))}} \rightarrow 0 \ \text{as} \ \Vert h_2 \Vert_{L^2(0,T;\mathbb{L}^2(\Omega))} \rightarrow 0. \]
Using the bilinear property of $ \widetilde{B} $, we proceed with
\begin{align*}
	&\widetilde{B}'\left[ S(\hat{h}+h_2), S(\hat{h}+h_2)+\alpha^2 A \, S(\hat{h}+h_2) \right] \tilde{u}\\
	& = \widetilde{B}\left[ \tilde{u},S(\hat{h}+h_2), S(\hat{h}+h_2)+\alpha^2 A \, S(\hat{h}+h_2) \right] +\widetilde{B}\left [S(\hat{h}+h_2), \tilde{u}+\alpha^2 A \tilde{u} \right ]\\
	&= \widetilde{B}\left[ \tilde{u}, \hat{u} + u_2 + r_1(h_2) + \alpha^2 A(\hat{u}+u_2+r_1(h_2)) \right] + \widetilde{B}\left[\hat{u}+u_2+r_1(h_2),\tilde{u}+\alpha^2 A \tilde{u} \right]\\
	&= \widetilde{B}'(\hat{u},\hat{u}+\alpha^2 A \hat{u}) \tilde{u} + \widetilde{B}'(u_2, u_2 + \alpha^2 Au_2)\tilde{u} + \widetilde{B}\left[r_1(h_2), r_1(h_2)+ \alpha^2 A r_1(h_2)\right]\tilde{u},
\end{align*}
where we employed $ u_2 = S'(\hat{h})h_2 $.
Thus, the difference $ \tilde{u} $ satisfies
	\begin{equation*}
	\begin{cases}
		\partial_{t} (\tilde{u} + \alpha^{2} A\tilde{u}) + \nu A(\tilde{u} + \alpha^{2} A\tilde{u}) + \widetilde{B}'(\hat{u},\hat{u}+\alpha^2 A \hat{u}) \tilde{u} \\
		\hspace{20.5pt} = h_1 - \widetilde{B}'(u_2, u_2 + \alpha^2 Au_2)\tilde{u} - \widetilde{B}(r_1(h_2), r_1(h_2)+ \alpha^2 A r_1(h_2))\tilde{u},\\
		\tilde{u}(0)=0.
	\end{cases}
\end{equation*}
Consequently, we obtain  \begin{equation} \label{S'tildeu}
	 \tilde{u} = S'(\hat{h}) \left[h_1 -\widetilde{B}'(u_2, u_2 + \alpha^2 Au_2)\tilde{u} - \widetilde{B}(r_1(h_2), r_1(h_2)+ \alpha^2 A r_1(h_2))\tilde{u} \right]. \end{equation}
On the other hand, let us estimate $ \widetilde{B}'(\hat{u},\hat{u}+\alpha^2 A\hat{u})u $ in the following way
\begin{equation*}
	\widetilde{B}'(\hat{u},\hat{u}+\alpha^2 A\hat{u})u = \widetilde{B}(u,\hat{u}+\alpha^2 A\hat{u}) + \widetilde{B}(\hat{u}, u+\alpha^2 Au),
\end{equation*}
and with $ v \in L^2(0,T;V) $
\begin{align*}
	|\widetilde{b}&(u, \hat{u}+\alpha^2 A\hat{u},v)| + |\widetilde{b}(\hat{u}, u+\alpha^2 Au, v)| \\
	& \leq |b(u,\hat{u}+\alpha^2 A \hat{u},v)| + |b(v,u,\hat{u}+\alpha^2 A \hat{u})|\\ 
	& \hspace{40.5pt}+ |b(\hat{u}, u+\alpha^2 Au, v)| + |b(v, \hat{u}, u+\alpha^2 Au)|\\
	&\leq C\Vert u \Vert |\hat{u}+\alpha^2 A \hat{u}| |Av| + C \Vert v \Vert |\hat{u}+\alpha^2 A\hat{u}| |Au| \\
	& \hspace{45.5pt}+ C \Vert \hat{u} \Vert |u+\alpha^2 Au| |Av| + C \Vert v \Vert |u+\alpha^2 Au| |A\hat{u}|\\
	&\leq C \Vert u \Vert \Vert \hat{u} \Vert \Vert v \Vert.
\end{align*}
This provides us with the estimation of the last term in \eqref{S'tildeu}
\begin{align*}
	& \left \Vert S'(\hat{h}) \left[ \widetilde{B}'(r_1(h_2), r_1(h_2)+\alpha^2 A r_1(h_2))\tilde{u}\right] \right \Vert_{W(0,T;H^3,V)}\\
	& \leq C \Vert r_1(h_2) \Vert_{W(0,T;H^3,V)} \Vert \tilde{u} \Vert_{W(0,T;H^3,V)} \\
	& \leq C \Vert r_1(h_2) \Vert_{W(0,T;H^3,V)} \Vert h_1 \Vert_{L^2(0,T;\mathbb{L}^2(\Omega))}.
\end{align*}
Next, we will investigate the second term of \eqref{S'tildeu}. By using again \eqref{S'tildeu} and the transformation $ \widetilde{B}'(u_2, u_2 + \alpha^2 A u_2) = \widetilde{B}''(\hat{u}, \hat{u} + \alpha^2 A \hat{u})(u_1, u_2) $, we can substitute $ \widetilde{B}'' $ instead of $ \widetilde{B}' $ and get
\begin{align*}
	\widetilde{B}'(u_2, u_2 + \alpha^2 &A u_2)\,\tilde{u}
	= \widetilde{B}''(\hat{u}, \hat{u}+\alpha^2 A \hat{u})(u_1, u_2) \\
	&- \widetilde{B}'(u_2, u_2 + \alpha^2 A u_2)\left\{ S'(\hat{h}) \left[ \widetilde{B}'(u_2, u_2 + \alpha^2 A u_2)\tilde{u} \right] \right\}\\
	&- \widetilde{B}'(u_2, u_2 + \alpha^2 A u_2)\left\{ S'(\hat{h}) \left[ \widetilde{B}'(r_1(h_2), r_1(h_2) + \alpha^2 A r_1(h_2))\tilde{u} \right] \right\}.
\end{align*}
Let us denote the second and third addend by $ r_2 $, e.g. 
\[ r_2 = \widetilde{B}'(u_2, u_2 + \alpha^2 A u_2)\left\{ S'(\hat{h}) \left[ \widetilde{B}'(u_2, u_2 + \alpha^2 A u_2)\tilde{u} + \widetilde{B}'(r_1(h_2), r_1(h_2) + \alpha^2 A r_1(h_2))\tilde{u} \right] \right\}. \]
It can be estimated as
\begin{align*}
	\Vert r_2 \Vert_{L^2(0,T;V')} &\leq C \Vert u_2 \Vert_{W(0,T;H^3,V)} \{ C_1 \Vert u_2 \Vert_{W(0,T;H^3,V)} \Vert \tilde{u} \Vert_{W(0,T;H^3,V)} \\
	&\hspace{15pt} + C_2 \Vert r_1(h_2) \Vert_{W(0,T;H^3,V)} \Vert \tilde{u} \Vert_{W(0,T;H^3,V)} \}\\
	&\leq C\Vert h_1 \Vert_{L^2(0,T;\mathbb{L}^2(\Omega))} \Vert h_2 \Vert_{L^2(0,T;\mathbb{L}^2(\Omega))} \Vert h_2 \Vert_{L^2(0,T;\mathbb{L}^2(\Omega))}\\
	&\hspace{15pt} +C\Vert h_1 \Vert_{L^2(0,T;\mathbb{L}^2(\Omega))} \Vert h_2 \Vert_{L^2(0,T;\mathbb{L}^2(\Omega))}\Vert r_1(h_2) \Vert_{W(0,T;H^3,V)}.
\end{align*}
Thus, we can write $ \tilde{u} $ as
\begin{align*}
	\tilde{u} = u_1 &- S'(\hat{h}) \left[\widetilde{B}''(\hat{u}, \hat{u} + \alpha^2 A \hat{u})(u_1,u_2) - r_2  \right]\\
	& - S'(\hat{h}) \, \widetilde{B}'(r_1(h_2), r_1(h_2)+\alpha^2 A r_1(h_2)) \, \tilde{u}.
\end{align*}
Hence, we find for the difference 
\begin{align*}
	\tilde{u} - u_1= - &S'(\hat{h}) \left[\widetilde{B}''(\hat{u}, \hat{u} + \alpha^2 A \hat{u})(u_1,u_2) - r_2  \right]\\
	& - S'(\hat{h}) \, \widetilde{B}'(r_1(h_2), r_1(h_2)+\alpha^2 A r_1(h_2)) \, \tilde{u}.
\end{align*}
Then we get for the remainder term
\begin{align*}
	&\frac{\left \Vert S'(\hat{h}+h_2)h_1 - S'(\hat{h})h_1- \left[-S'(\hat{h})\,\widetilde{B}''(\hat{u},  \hat{u}+\alpha^2 A\hat{u})(u_1,u_2)\right] \right \Vert_{W(0,T;H^3,V)}}{\Vert h_2 \Vert_{L^2(0,T;\mathbb{L}^2(\Omega))}}\\
	&= \frac{\left \Vert \tilde{u} - u_1- \left[-S'(\hat{h})\,\widetilde{B}''(\hat{u},  \hat{u}+\alpha^2 A\hat{u})(u_1,u_2)\right] \right \Vert_{W(0,T;H^3,V)}}{\Vert h_2 \Vert_{L^2(0,T;\mathbb{L}^2(\Omega))}}\\
	&\leq \frac{\Vert S'(\hat{h})r_2 \Vert_{W(0,T;H^3,V)} + \Vert S'(\hat{h})\widetilde{B}'(r_1(h_2), r_1(h_2)+\alpha^2 A r_1(h_2)) \Vert_{W(0,T;H^3,V)}}{\Vert h_2 \Vert_{L^2(0,T;\mathbb{L}^2(\Omega))}}\\
	&\leq \frac{C \Vert h_1 \Vert_{L^2(0,T;\mathbb{L}^2(\Omega))} \Vert h_2 \Vert_{L^2(0,T;\mathbb{L}^2(\Omega))} \left[\Vert h_2 \Vert_{L^2(0,T;\mathbb{L}^2(\Omega))} + \Vert r_1(h_2) \Vert_{W(0,T;H^3,V)}\right]}{\Vert h_2 \Vert_{L^2(0,T;\mathbb{L}^2(\Omega))}}\\
	&\hspace{20pt} + \frac{\Vert r_1(h_2) \Vert_{W(0,T;H^3,V)} \Vert h_1 \Vert_{L^2(0,T;\mathbb{L}^2(\Omega))}}{\Vert h_2 \Vert_{L^2(0,T;\mathbb{L}^2(\Omega))}},
\end{align*}
which converges to $ 0 $ as $ \Vert h_2 \Vert_{L^2(0,T;\mathbb{L}^2(\Omega))} \rightarrow 0 $. This infers that the second derivative of $ S $ is given as
\[ S''(\hat{h})(h_1, h_2) = -S'(\hat{h})\widetilde{B}''(\hat{u}, \hat{u}+\alpha^2 A \hat{u})(u_1, u_2), \]
or equivalently as the weak solution of \eqref{2nd_S}.\\
It remains to prove the continuity of $ S''(h) $ with respect to $ h $. We shall decompose $ S''(\hat{h})(h_1, h_2) $ as
\begin{align*}
	S''(\hat{h})(h_1, h_2) &= -S'(\hat{h}) \, \widetilde{B}''(\hat{u}, \hat{u}+ \alpha^2 A \hat{u})(u_1, u_2)\\
	&= -S'(\hat{h})\, \left[ \widetilde{B}(u_1, u_2 + \alpha^2 A u_2) + \widetilde{B}(u_2, u_1 + \alpha^2 A u_1) \right].
\end{align*}
Here, the solution mappings $ S(\hat{h}) $ and $ S'(\hat{h}) $ depend continuously on the parameter $ \hat{h} $, see Theorem \ref{1stS}. Altogether the mapping $ \hat{h} \mapsto S''(\hat{h})[h_1, h_2] $ is continuous and $ S(h) $ belongs to class $ C^2 $.
\end{proof}

\subsection{Assumptions on the objective functions and some obtained results}
For the convenience of the reader, we shall recall some assumptions imposed on the objective function $ J $ and the obtained results in \cite{AGN2023}.

We call a couple $ (u,h) $ of state and control admissible if it satisfies the constraints \eqref{state equation} and \eqref{pointwise_constraint} of the optimal control problem. Moreover, we denote by $$\mathcal{A} = \{(u,h) \in \mathcal{H}\times W(0,T;H^3,V): u  \text{ is the unique weak solution to } \eqref{state equation} \} $$ the set of admissible pairs.\\
Sometimes, it is convenient to work with the reduced objective functional $ J(h) $ that is defined by
\[ J(h)=J(S(h),h), \]
where $ S: h\mapsto u $ is the (nonlinear) solution operator associated to \eqref{state equation}. 

Before confirming the existence of optimal solutions to the problem $(\textbf{P})$ and proving the first-order necessary optimality conditions, we need the following assumptions.

$(\textbf{A}_1)\mbox{ } F$ and $L$ are Carathéodory functions. Besides, we assume that for every $(x,t,u)\in Q\times \mathbb{R}^3$ and $ x \in \Omega$ respectively, $L(x,t,u,\cdot)$ and also $F(x,\cdot)$ are convex functions.

$(\textbf{A}_2)$  The functions $ F $ and $ L $ satisfy the local Lipschitz conditions, that is, for any $M>0$ there exists functions  $l(x,t) \in L^{\infty}(Q,\mathbb{R})$ and $ l'(x,t) \in L^{\infty}(\Omega,\mathbb{R})$ such that the following properties holds
\begin{align}
	|L(x,t,u_1,h_1) - L(x,t,u_2,h_2)| &\leq l(x,t).(|u_1-u_2|+|h_1-h_2|), \label{Liptchiz_L}\\
	|F(x, u_1(x,T)) - F(x,u_2(x,T))|&\leq l'(x,t).(|u_1(x,T)-u_2(x,T)|) \label{Lipchitz_F}
\end{align}
for a.e $ (x,t) \in Q $ and $ |u_1|, |h_1|,|u_2|, |h_2| \leq M $.

$(\textbf{A}_3)$  There exist $\psi \in L^1(Q),\phi\in L^1(\Omega) $ and $C> 0$ such that
\[ L(x, t, u,h) \geq C|h|^2+\psi (x,t), \ \forall (x,t,u,h) \in Q \times \mathbb{R}^3 \times \mathbb{R}^3. \]
\[ F(x, u)\geq \phi (x), \ \forall (x,u) \in \Omega \times \mathbb{R}^3 . \]

$(\textbf{B}_1)$ Beside the Carathéodory condition, we assume that  $F(x, \cdot)$ and $ L(x,t,\cdot,h(x,t))$ are differentiable respectively around $\hat{u}(x,T),\ \hat{u}(x,t)$ for a.e $(x,t)\in Q$ and $h\in\mathcal{H}$. 

$(\textbf{B}_2)$ There exist functions $l_1(x,t) \in L^{\infty}(Q,\mathbb{R})$ and $l_2(x,t) \in L^{\infty}(\Omega,\mathbb{R})$ such that
\begin{align}
	|L_u(x,t,u_1,h_1) - L_u(x,t,u_2,h_2)|&\leq l_1(x,t).(|u_1-u_2|+|h_1-h_2|) \label{Liptchiz_Lu}\\
	|F_u(x, u_1(x,T)) - F_u(x,u_2(x,T))|&\leq l_2(x,t).(|u_1(x,T)-u_2(x,T)|) \label{Lipchitz_Fu}
\end{align}
for a.e $(x,t)\in Q$ and for all $(u_1,h_1), (u_2,h_2)$ in a neighborhood of $(\hat{u},\hat{h})$. Here, $ F_u $ and $ L_u $ denote the derivatives of $ F $ and $ L $ with respect to $ u $, respectively.

\begin{theorem} \cite{AGN2023}\label{exist sol}
	There exists at least one optimal control solution $\hat{h}\in \mathcal{H}$ to $(\textbf{P})$ with associated state $\hat{u}\in L^\infty(0,T;\mathbb{H}_{\sigma}^2(\Omega))\cap W(0,T;H^3,V)$.
\end{theorem} 
The first-order derivatives of $ J $ with respect to $ u$ and $ h $ in directions $ w \in W(0,T;H^3,V) $ and $ k \in \mathbb{L}^2(Q) $ respectively are
\begin{subequations}
	\begin{align}
		& J_u(u,h) w  = \int_Q L_u(x,t,u,h)w \, dxdt + \int_\Omega F_u(x,u(x,T)) w \, dx, \label{1st u} \\
		&J_h(u,h)k = \int_Q L_h(x,t,u,h)k \, dxdt.  \label{1st h}
	\end{align}
\end{subequations}

Before stating the first-order condition in this section, we recall some definitions of Convex Analysis. For a convex subset $ U $ of a Hilbert space $ H $ and an element $ h \in H $, we denote by $ \mathcal{N}_U (h) $ and $ \mathcal{T}_U(h) $ the normal cone and the polar cone of tangents of $ H $ at the point $ h \in H $, respectively, i.e.
\[ \mathcal{N}_U(h) = \{z\in H : (z, v - h) \leq 0 \ \forall v \in U\}, \]
\[ \mathcal{T}_U(h) = \{z\in H : (z, v) \leq 0 \ \forall v \in \mathcal{N}_U(h) \}. \]
An element $ \omega \in H $ is called a feasible direction at $ h \in U $ if there exists $ \delta > 0 $ such that $ h + \epsilon \omega \in U$ holds for all $ \epsilon \in (0, \delta) $.
The cone of feasible directions at $ h\in U $, denoted by $ \mathcal{F}_U(h) $, has been proved (see e.g. \cite{Aubin}) to possess the following property
\begin{equation} \label{closure_feasiblecone}
	\overline{\mathcal{F}_U(h)} =  \mathcal{T}_U(h). 
\end{equation}
We apply the mentioned definitions for the case $H = \mathbb{L}^2(Q) $, $ U = \mathcal{H} $ and $ h = \hat{h} $.

\begin{theorem} \label{FNC}
	Let $ \hat{h} $ be locally optimal in $ \mathbb{L}^2(Q) $ with associated state $ \hat{u} $. 
	Then there exists a unique weak solution  $\hat{\lambda} \in  L^\infty(0,T;\mathbb{H}_{\sigma}^2(\Omega)) \cap L^2(0,T;\mathbb{H}_{\sigma}^3(\Omega))$ with $ \partial_t \hat {\lambda} \in L^2(0,T;V) $ of the adjoint equations
	\begin{equation}\label{adjoint eq}
		\begin{cases}
			-\partial_t(\hat{\lambda}+\alpha^2A\hat{\lambda}) +\nu A(\hat{\lambda}+\alpha^2A\hat{\lambda})+\hat{B}(\hat{u},\hat{\lambda}) = L_u(.,t,\hat{u}(.,t),\hat{h}(.,t))\\
			\hfill \text{ in } V' \text{ for a.e. } t\in [0,T],\\
			A\hat{\lambda} = 0,\hfill x\in\partial \Omega, t\in [0,T],\\
			\hat{\lambda}(T)+\alpha^2 A\hat{\lambda}(T) = F_u(.,\hat{u}(.,T)) \text{ in } V',
		\end{cases}
	\end{equation}
	where $\hat{B}(\hat{u},\hat{\lambda})\in V'$ is defined as follows
	$$\langle \hat{B}(\hat{u},\hat{\lambda}),w\rangle_{V',V}:=\tilde{b}(\hat{u},w+\alpha^2Aw,\hat{\lambda})+\tilde{b}(w,\hat{u}+\alpha^2A\hat{u},\hat{\lambda}).$$
	The other terms in the first and last equations of \eqref{adjoint eq} take their dualities with functions in $V$ by using the inner product in $\mathbb{L}^2(\Omega)$.\\
 Moreover, we have
	\begin{equation} \label{origin_variational_ineq}
		J_h(\hat{u},\hat{h}) \widetilde{h}   \geq 0 \ \forall \widetilde{h} \in  \mathcal{T}_{\mathcal{H}}(\hat{h}).  
	\end{equation}
	As a special case, the inequality
	\begin{equation} \label{variational ineq}
		J_h(\hat{u},\hat{h}) (h - \hat{h})   \geq 0, \ \forall h \in \mathcal{H}, 
	\end{equation}
	is satisfied.
\end{theorem}

\begin{remark}
	\rm{In fact, inequalities \eqref{origin_variational_ineq} and \eqref{variational ineq} were established in \cite{AGN2023} for the Lagrange functional 
		 \[ \mathcal{L}: W(0,T;H^3,V)\times \mathbb{L}^2(Q) \times L^2(0,T;\mathbb{H}_{\sigma}^3(\Omega)) \to \mathbb{R} \]
		 defined by
		 \begin{align*} 
		 	&\mathcal{L}(u,h,\lambda)=  J(u,h) \\
		 	&- \int_0^T\left \langle \partial_{t} (u + \alpha^{2} Au) + A(u + \alpha^{2} Au) + \widetilde{B}(u,u + \alpha^{2} Au) - h, \lambda \right \rangle _{V', V}dt.
	 \end{align*}	
 However, by arguments therein, we can prove a similar estimations for $ J_h(u, h) $.
 }
\end{remark}

\section{SECOND-ORDER SUFFICENT OPTIMALITY CONDITIONS}
Let $ \hat{h} $ be an $ L^2 $-locally optimal solution of our optimal control problem and we assume that 
the functions $ F $ and $ L $ are twice continuously differentiable respectively around the optimal points $ (\hat{u}, \hat{h})$.
We will prove that the second derivative of the objective function $ J $ is positive on the cone of critical directions defined below.
\begin{definition}(Cone of critical directions)\\
	Using the reduced objective functional, we denote by $ C_0({\hat{h}}) $ the cone of critical directions, meaning that
	\[ C_0({\hat{h}})=\{\widetilde{h}\in \mathcal{T}_{\mathcal{H}}(\hat{h}): J'(\hat{h})\widetilde{h}=0\}. \]
\end{definition}
Moreover, we establish the second-order derivative of the objective functional $ J $ with respect to $ v=(u,h) $ in directions $ (w_1,k_1) $ and $ (w_2,k_2) \in W(0,T;H^3,V) \times \mathbb{L}^2(Q)$.
\begin{align*}
	& J_{vv}(u,h) [(w_1,k_1),(w_2,k_2)]  = \int_Q L_{vv}(x,t,u,h)[(w_1,k_1),(w_2,k_2)] \, dxdt  \\
	&+ \int_\Omega F_{uu}(x,u(x,T)) w_1(x,T)w_2(x,T) \, dx. 
\end{align*}

The main result of this section is the following theorem.
\begin{theorem}\label{SSC}
	Let $ \hat{v} = (\hat{u}, \hat{h}) $ be an admissible pair and suppose that $ \hat{v} $ satisfies the first-order necessary optimality conditions. Moreover, we assume that the pair $ \hat{v}=(\hat{u},\hat{h}) $ satisfies the following assumption, known as the second-order sufficient condition (SSC). It holds
	\begin{equation}\label{SSC_full}
		J_{vv}(\hat{u},\hat{h}) [(z,\widetilde{h})^2] > 0, 
	\end{equation}
or equivalently with the reduced objective function
\begin{equation}\label{SSC}
	J''(\hat{h}) \widetilde{h}^2 > 0
\end{equation}
	for all $ \widetilde{h} \in C_0({\hat{h}}) $, where $ z $ is the unique solution of the linearized equation \eqref{linearized eq} with the right-hand side $\widetilde{h}  $.
	Then there exist $ \epsilon >0 $ and $ \rho>0 $ such that 
	\begin{equation} J(h)\geq J(\hat{h})+\epsilon \, \| h_1 - h_2\|_{\mathbb{L}^2(Q)}^2 \label{growth_ineq}
	\end{equation} 
	holds for all $ h\in \mathcal{H} $ with $ \| h_1 - h_2\|_{\mathbb{L}^2(Q)}^2 \leq \rho $, which implies that $ \hat{h} $ is a locally optimal control with associated state $ \hat{u} $ .
\end{theorem} 
\begin{proof}
	Let us suppose that the first-order necessary and second-order sufficient conditions are satisfied, whereas \eqref{growth_ineq} does not hold. Then for all $ \epsilon > 0 $ and $ \rho >0 $ there exists $ h_{\epsilon, \rho } \in \mathcal{H}  $ with $ \| h_{\epsilon, \rho } - \hat{h} \|_{\mathbb{L}^2(Q)}  \leq \rho$ and 
	\[ J(h_{\epsilon, \rho }) < J(\hat{h}) +\epsilon \, \| h_{\epsilon, \rho } - \hat{h} \|_{\mathbb{L}^2(Q)}^2, \]
	where $ u_{\epsilon, \rho } $ is the corresponding state of $ h_{\epsilon, \rho } $.\\
	Hence, for any $ k \in \mathbb{Z}^{+} $, we have admissible pairs $ (u_k,h_k) $ such that 
	\begin{equation} \label{contradicts_growth_ineq}
		J(h_k) < J( \hat{h}) + \frac{1}{k} \, \| h_{k } - \hat{h} \|_{\mathbb{L}^2(Q)}^2   
	\end{equation}
	and $ \| h_{k } - \hat{h} \|_{\mathbb{L}^2(Q)} < \frac{1}{k} $.
	By construction, we have $ h_k \rightarrow \hat{h} $ in $ \mathbb{L}^2(Q) $ as $ k \rightarrow \infty $.
	Thus, we can write $ h_k = \hat{h} + t_k q_k $, where $ t_k>0  $ and $ q_k \in \mathcal{F}_{\mathcal{H}}(\hat{h}), \|q_k\|_{\mathbb{L}^2(Q)} =1  $ and $ t_k \rightarrow 0 $ as $ k \rightarrow \infty $.\\
	With application of arguments in Lemma \ref{LipchitzS},	and the assumption $ \|q_k\|_{\mathbb{L}^2(Q)} =1 $, $ \{z_k\} $ is bounded in $ W(0,T;H^3,V) $, where $ z_k $ is the unique solution of linearized equations with the right-hand side $ q_k $.
	Hence, we can extract subsequences denoted again by $ \{q_k\} $ and $ \{z_k\} $ converging weakly to some $ \widetilde{h} \ \text{in} \ \mathbb{L}^2(Q) $ and $ \widetilde{z} \ \text{in} \ W(0,T;H^3,V) $. Since $ W(0,T;H^3,V) $ is compactly embedded in $ \mathbb{L}^2(Q) $, we obtain $ z_k \rightarrow \widetilde{z} $ in $ \mathbb{L}^2(Q) $ as $ k \rightarrow \infty $. It follows that $ \widetilde{z}  $ is the unique solution of the linearized equations with the right-hand side of the first equation $ \widetilde{h} $. It remains to show that $ \widetilde{h} \in C_0(\hat{h}) $ and $ J''(\hat{h})\widetilde{h}^2 \leq 0 $, which contradicts \eqref{SSC} and so we get the claim.\\
	Since $ \mathcal{T}_{\mathcal{H}}(\hat{h}) $ is convex and closed, we have $ \widetilde{h} \in \mathcal{T}_{\mathcal{H}}(\hat{h}) $ and thus $ J'(\hat{h})\widetilde{h} \geq 0 $ . 
	Meanwhile, since $ J $ is twice continuously differentiable around optimal points
	\begin{align*}
		J'(\hat{h})\widetilde{h} &= \lim _{k\rightarrow\infty} \dfrac{J(\hat{h}+t_kq_k) - J(\hat{h})}{t_k} = \lim _{k\rightarrow\infty}\dfrac{J(h_k) - J(\hat{h})}{t_k}\\
		&\leq\lim _{k\rightarrow \infty}\frac{1}{k} \, \| h_k - \hat{h} \|_{\mathbb{L}^2(Q)} \leq \lim _{k\rightarrow \infty}\frac{1}{k^2} = 0,
	\end{align*}
which implies that $ h \in C_0(\hat{h}) $. Finally, we show that $ J''(\hat{h})\widetilde{h}^2 \leq 0 $. Using Taylor's formula we obtain
	\[ J(\hat{h}+t_kq_k) = J(\hat{h})+t_k J'(\hat{h})q_k + \frac{t_k^2}{2} J''(\hat{h}+\theta_k t_k q_k)q_k^2,   \]
	with $ \theta_k \in (0,1) $.
	Thus
		\[ \frac{t_k^2}{2}J''(\hat{h}+\theta_k t_k q_k)q_k^2 \leq J(h_k) - J(\hat{h})
		< \frac{1}{k} \, \| h_k - \hat{h} \|_{\mathbb{L}^2(Q)}  < \frac{t_k^2}{k}. \]
	Therefore, it holds that
	\begin{equation} J''(\hat{h}+\theta_k t_k q_k)q_k^2 < \frac{2}{k}. \label{J''} \end{equation}
	Using again the fact that $ J $ is twice continuously differentiable around optimal points $ (\hat{u},\hat{h}) $, and $\hat{h}+\theta_k t_k q_k  $ strongly converges to $ \hat{h} $ as $ k \rightarrow \infty $, we get
	\begin{align}  
	  & \left|  \left[ J''(\hat{h}) - J''(\hat{h}+\theta_k t_k q_k) \right]  q_k^2 \right| \nonumber \\
	  &\leq \sup \left\{ \left| \left[ J''(\hat{h}) - J''(\hat{h}+\theta_k t_k q_k) \right] uv \right|: \Vert u \Vert = \Vert v \Vert =1 \right\} \nonumber \\
	  &= \left \Vert J''(\hat{h}) - J''(\hat{h}+\theta_k t_k q_k) \right \Vert \rightarrow 0 \ \text{as} \ k \rightarrow \infty. \label{supJ}
	\end{align}
	 Collecting \eqref{J''} and \eqref{supJ} we obtain $ J''(\hat{h})\widetilde{h}^2 \leq 0 $, which infers the desired contradiction.
\end{proof}
\begin{remark}
	{\rm \textbf{1.} In case of box constraints, where $ \mathcal{H} $ is defined by \[ \mathcal{H} = \{h\in \mathbb{L}^2(Q): h_{a,i}(x,t) \leq h_i (x,t)\leq h_{b,i}(x,t) \ \text{a.e on} \ Q, i =1,2 \}, \] 
	by constructing the cone of $ L^{\infty}$-functions $ \widetilde{C}_0(\hat{h}) = C_0(\hat{h}) \cap (L^{\infty}(Q))^3 $ and approximating $ \widetilde{C}_0(\hat{h}) $	by a family of cones $ \widetilde{C}_{\sigma}(\hat{h})  $ as in \cite{Wachsmuth2006}, we can prove the second-order necessary optimality conditions as the following
	\[ J_{vv}(\hat{u},\hat{h})[(z,\widetilde{h})^2] > 0 \ \text{for every} \ \widetilde{h} \in C_0(\hat{h}), \]
	and $ z$ is the solution of the linearized equations with the right-hand side $ 
	\widetilde{h} $.\\
	
	\textbf{2.} Considering the \textit{quadratic form} of the objective functional $ J $ mentioned in \cite{AS2019}
	\begin{align*}
	J(u,h) =\frac{\alpha_Q}{2} & \int_Q  |u(x,t)-u_d(x,t)|^2 dxdt +\frac{\alpha_T}{2} \int_{\Omega} |u(x,T)-u_T(x)|^2 \ dx \\&+ \frac{\gamma}{2} \int_Q|h(x,t)|^2 dxdt,  
	\end{align*}
	the cone of critical directions is rewritten as
	\[ C_0(\hat{h}) = \left \{\widetilde{h} \in \mathcal{T}_{\mathcal{H}}: \int_Q (\hat{\lambda} + \gamma \hat{h})\widetilde{h} \ dxdt = 0 \right \}. \]
	Thus, we get the $ (SSC) $
	\begin{align*} \frac{\alpha_T}{2} & \int_\Omega |z(x,T)|^2 dx + \frac{\alpha_Q}{2} \int_Q |z(x,t)|^2 dxdt + \frac{\gamma}{2} \int_Q |\widetilde{h}|^2 dxdt \\& - \int_{0}^{T} \widetilde{b}(z(t),z(t)+\alpha^2 Az(t), \hat{\lambda}(t)) dt >0, 
	\end{align*}
	where $ z$ is the solution of the linearized equations with the right-hand side $ 
	\widetilde{h} $. This is the exactly the second-order sufficient optimality conditions obtained in \cite{AS2019}.
}
\end{remark}
\section{STABILITY OF THE OPTIMAL CONTROL PROBLEMS WITH RESPECT TO INITIAL DATA }
In this section,  following the general lines of the recent approach introduced by Casas and Tr\"{o}ltzsch in \cite{CT2022},  we will study the stability of selected local solutions of the optimal control problem with respect to a perturbation of the initial data $ u_0 $. Let $ \hat{h} $ be a fixed local solution of problem $ (\textbf{P}) $, we estimate the distance from $ \hat{h} $ to an associated local minimizer of the perturbed problem with the initial function $ u_0 + \phi_\varepsilon$, where $ \left \| \phi_\varepsilon \right \|  $ is small enough.

First, we reformulate the state equations \eqref{state equation} by using the operators defined in Section 2 as follows
\begin{equation}\label{VCHE_original}
	\begin{cases}
		\partial_{t}(u + \alpha^2 Au) + \nu A(u +\alpha^2 Au) + \widetilde{B}(u,u+\alpha^2 Au) = h \ \text{in} \ L^2(0,T;V'),\\
		u(0)=u_0 \ \text{in} \ D(A).
	\end{cases}
\end{equation}
Now, we consider perturbations in the initial condition of \eqref{state equation} leading to a family of perturbed optimal control problem $ (\textbf{P}_\varepsilon) $. Let $ \{ \phi _\varepsilon\}_{\varepsilon>0} \subset D(A)$ be a family of functions satisfying 
\begin{align}
	&\exists M_\phi < \infty \ \text{such that} \ \left \| \phi_\varepsilon \right \|
	\leq M_\phi \ \forall \varepsilon >0, \label{1st_perturbation}\\
	&\lim_{\varepsilon \rightarrow 0}\left \| \phi_\varepsilon \right \| = 0. \label{2nd_perturbation}
\end{align}
Remark that we need $ \phi_\varepsilon$ in $D(A) $ to ensure that the perturbed equations has a unique weak solution $ u $ in $ (L^\infty (Q))^3$. However, the requirement $ \underset{\varepsilon \rightarrow 0}\lim \left \| \phi_\varepsilon \right \|_{D(A)} = 0 $ is too strong, so the selection of the $ V $-norm is to obtain more practical results. 

We associate with this family the state equations
\begin{equation*}
	\begin{cases}
		\partial_{t}\left ( u-\alpha ^{2} \Delta u \right )+\nu  \left (Au-\alpha ^{2} \Delta Au  \right ) + \nabla p\\
		\hspace{70pt}= u \times \left ( \nabla\times\left ( u-\alpha ^{2}\Delta u \right ) \right ) + h,  &(x,t) \in (\Omega \times (0,+\infty)),\\
		\nabla \cdot u = 0, &(x,t) \in (\Omega \times (0,+\infty)), \\
		u = Au = 0,   &(x,t) \in (\partial \Omega \times (0,+\infty)),  \\
		u(x,0) = u_0(x) + \phi_\varepsilon,  &\ x\in\Omega,
	\end{cases}
\end{equation*}
or equivalently,
\begin{equation}\label{VCHE_perturbed}
	\begin{cases}
		\partial_{t}(u + \alpha^2 Au) + \nu A(u +\alpha^2 Au) + \widetilde{B}(u,u+\alpha^2 Au) = h \ \text{in} \ L^2(0,T;V'),\\
		u(0)=u_0 +\phi_\varepsilon \ \text{in} \ D(A).
	\end{cases}
\end{equation}
For given $ \varepsilon $ and $ u $, the solution of this equation will be denoted by $ u_h^\varepsilon $. Then, we consider the pertubed optimal control problems
\begin{equation} \label{objective_function_perturbed} \tag{$ \textbf{P}_\varepsilon $}
	 \min_{h\in \mathcal{H}} J_\varepsilon(h) = \int_Q L(x,t,u_h^\varepsilon(x,t),h(x,t))dxdt+\int_{\Omega}F(x,u_h^\varepsilon(x,T))dx.
\end{equation}
Similarly to problem $ (\textbf{P}) $, every problem $ (\textbf{P}_\varepsilon) $ has at least one global solution $ h_\varepsilon $. \\
The next two theorems analyze the relation between the solutions of $ (\textbf{P}) $ and $ (\textbf{P}_\varepsilon) $.
\begin{theorem} \label{1st_theorem_stability}
	Let $ \{h_\varepsilon\}_{\varepsilon >0} $ be a sequence of global solutions of problems $ (\textbf{P}_\varepsilon) $. Any control $ \hat{h} $ that is the weak limit in $ \mathbb{L}^2(Q) $ of a sequence $ \{ h_{\varepsilon_k}\}_{k=1}^\infty$ with $ \varepsilon_k \rightarrow 0$ as $ k \rightarrow \infty $ is a global minimizer of $ (\textbf{P}) $. Moreover, the convergence is strong in $ \mathbb{L}^2(Q) $.
\end{theorem}
\begin{proof}
	First, we shall prove the uniform boundedness of $ h_\varepsilon $ in $ \mathbb{L}^2(Q) $. 
	Utilizing assumption $(\textbf{A}_3)$ and the fact that $ h_\varepsilon \in \mathcal{H}$, we deduce
	\begin{align*}
		C \Vert h_\varepsilon \Vert_{\mathbb{L}^2(Q)} ^2 + \int_Q \psi(x,t) \,dxdt &\leq \int_Q L(x,t,u,h) \, dxdt \\ &\leq J - \int_{\Omega} F(x,u(x,T))dx \leq J - \int_{\Omega}\phi(x)\,dx.
	\end{align*}
	Thus,
	\begin{equation} \label{h_bounded}
		C \Vert h_\varepsilon \Vert_{\mathbb{L}^2(Q)}^2 + \int_Q \psi(x,t)\,dxdt + \int_{\Omega}\phi(x)\,dx \leq J_\varepsilon(h_\varepsilon).
	\end{equation}
	Since $ h_\varepsilon $ is a global solution of problem $ (\textbf{P}_\varepsilon) $, one can find $ h_0 \in \mathcal{H} $ such that $ J_\varepsilon (h_\varepsilon) \leq J_\varepsilon(h_0)$.\\
	Considering $ J_\varepsilon(h_0) - J(h_0) $, we obtain the following estimates
	\begin{align*}
		\vert J_\varepsilon(h_0) - J(h_0) \vert &= \int_Q \vert L(x,t,u^\varepsilon_{h_0}(x,t), h_0(x,t)) - L(x,t,u_{h_0}(x,t),h_0(x,t))\vert \,dxdt 
		\\& \hspace{24pt}+ \int_{\Omega} \vert F(x,u^\varepsilon_{h_0}(x,T)) - F(x,u_{h_0}(x,T))\vert \, dx \\
		& \leq \int_Q C \vert u^\varepsilon_{h_0} - u_{h_0} \vert \, dxdt + \int_{\Omega} C \vert u^\varepsilon_{h_0}(x,T) - u_{h_0}(x,T) \vert \,dx \\
		& \leq \int_Q C \Vert \phi_\varepsilon \Vert \,dxdt + \int_{\Omega} C \Vert \phi_\varepsilon \Vert \,dx \leq C,
	\end{align*}
	where $ u^\varepsilon_{h_0} $ and $ u_{h_0} $ denotes the solutions of the perturbed and unperturbed equations with the right-hand side $ h_0 $, respectively. \\
	Furthermore, since $L: Q\times\mathbb{R}^3\times\mathbb{R}^3 \rightarrow \mathbb{R}$ and $F: \Omega \times \mathbb{R}^3 \rightarrow \mathbb{R}$, we get  \[ \vert J_\varepsilon(h_0) \vert \leq C + \vert J(h_0) \vert < K . \]
	Together with \eqref{h_bounded}, we obtain the uniform boundedness of $ h_\varepsilon $ in $ \mathbb{L}^2(Q) $, which implies the existence of the sequences $ \{ h_{\varepsilon_k}\}_{k=1}^{\infty} $ converging weakly to $ \hat{h} $ in $ \mathbb{L}^2(Q) $. \\
	We denote by $ u_{\varepsilon_k} $ and $ \hat{u} $ the associated states with $ h_{\varepsilon_k} $ and $ \hat{h} $, solutions of \eqref{VCHE_perturbed} and \eqref{VCHE_original}, respectively. Computing similarly as in Lemma \ref{LipchitzS}, 
we infer that $ u_{\varepsilon_k} $ converges weakly to $ \hat{u} $ in $ W(0,T;H^3,V) $ and strongly in $ \mathbb{L}^2(Q) $.
Using the optimality of $ h_{\varepsilon_k} $, together with assumption $(\textbf{A}_2)$ and arguing as [\cite{AGN2023}, Theorem 3.1], we deduce for every $ h \in \mathcal{H} $
\[ J(\hat{h}) \leq \liminf_{k \rightarrow \infty } J_{\varepsilon_k}(h_{\varepsilon_k}) \leq \limsup_{k\rightarrow \infty}J_{\varepsilon_k}(h_{\varepsilon_k}) = J(h).   \]
Since $ \hat{h} \in \mathcal{H} $, the above inequalities imply that $ \hat{h} $ is a global solution of $ (\textbf{P}) $. Then the strong convergence of $ h_{\varepsilon_k} $ to $ \hat{h} $ in $ \mathbb{L}^2(Q) $ follows from the convergence $ u_{\varepsilon_k} \rightarrow \hat{u} $ and assumptions $ (\textbf{A}_2)$, $ (\textbf{A}_3)$.
\end{proof}
Conversely, following arguments in \cite{CT2022}, we can prove the following theorem.
\begin{theorem} \label{2nd_theorem_stability}
	Let $ \hat{h} $ be a strict local minimizer of $ (\textbf{P}) $, meaning that $ J(\hat{h}) < J(h) $ whenever $ h \neq \hat{h} $. Then, there exists a set $ \{ h_\varepsilon\}_{\varepsilon >0} $ of local solutions of the perturbed problems $ (\textbf{P}_\varepsilon) $ such that $ h_\varepsilon \rightarrow \hat{h} $ strongly in $ \mathbb{L}^2(Q) $ when $ \varepsilon \rightarrow 0 $.
\end{theorem}
\begin{proof}
	Since $\hat{h}$ is a strict local solution of $ (\textbf{P}) $, there exists a closed $ \mathbb{L}^2(Q) $-ball $ B_\rho(\hat{h}) $ such that $ J(\hat{h}) < J(h) $ for every $ u \in \mathcal{H} \cap B_\rho(\hat{h}) \setminus \{\hat{h}\}$. We consider the control problems
	\begin{equation} \tag{$ \mathcal{P} $}
	\min_{h\in \mathcal{H} \cap B_\rho(\hat{h})} J(h)	
	\end{equation}
and \begin{equation}\tag{$ \mathcal{P}_\varepsilon $}
	\min_{h\in \mathcal{H} \cap B_\rho(\hat{h})} J_\varepsilon(h).
\end{equation}
	Obviously, $ \hat{h} $ is the unique solution of $( \mathcal{P} )$ and every problem $ (\mathcal{P}_\varepsilon) $ has at least one solution $ h_\varepsilon $. From Theorem \ref{1st_theorem_stability}, any control that is the weak limit in $ \mathbb{L}^2(Q) $ of the sequence $ \{ h_{\varepsilon_k}\}_{k=1}^\infty $ is a global minimizer of $ (\textbf{P})$. However, since $ \hat{h} $ is the unique solution of $( \mathcal{P} )$, the whole family $ \{ h_\varepsilon \}_{\varepsilon>0} $ converges to $ \hat{h} $ and this convergence is strong in $ \mathbb{L}^2(Q) $. This implies that there exists $ \varepsilon_0>0 $ such that $ \vert h_\varepsilon - \hat{h} \vert < \rho  $ for every $ \varepsilon < \varepsilon_0 $. Therefore, $ h_\varepsilon $ is a local solution of $ (\textbf{P}_\varepsilon) $ for every $ \varepsilon < \varepsilon_0 $.
\end{proof}

Let $ \hat{h} $ be the solution to problem $ (\textbf{P}) $ satisfying the second-order sufficient optimality condition \eqref{SSC}. Theorem \ref{2nd_theorem_stability} establishes the existence of a set $ \{ h_\varepsilon\}_{\varepsilon >0} $ of solutions to $ (\textbf{P}_\varepsilon) $ such that $ h_\varepsilon \rightarrow \hat{h} $ in $ \mathbb{L}^2(Q) $ as $ \epsilon \rightarrow 0 $. Finally, we estimate $ h_\varepsilon - \hat{h} $.
\begin{theorem} \label{3rd_theorem_stability}
	Let $ \hat{h} $ be a solution of the optimal control problem satisfying the second-order sufficient condition \eqref{SSC}
	\[ J''(\hat{h})\widetilde{h}^2 > 0 \ \text{with} \ \widetilde{h} \in C_0(\hat{h}).\]
	Then there exists $ \varepsilon_0 >0$ and a constant $ C  $ such that
	\begin{equation} \label{lipchitz_convergence}
		\Vert h_{\varepsilon} - \hat{h} \Vert_{\mathbb{L}^2(Q)} \leq C \Vert \phi_\varepsilon \Vert_{\mathbb{L}^2(\Omega)} \ \forall \varepsilon \in (0,\varepsilon_0).
	\end{equation}
\end{theorem}
Before proving the above theorem, we need some auxiliary lemmas.
\begin{lemma} \label{lemma1}
	Let $ u_{h_\varepsilon} $ and $ u^{h_\varepsilon} $ denote the solutions of the unperturbed equations \eqref{VCHE_original} and the perturbed equations \eqref{VCHE_perturbed} with the right-hand side $ h_\varepsilon $. Similarly, $ \lambda_{h_\varepsilon} $ and $ \lambda^{h_\varepsilon} $ are the corresponding adjoint states. Then there exists a constant $ C > 0 $ 
	such that
	\begin{equation}
		\Vert u_{h_\varepsilon} - u^{h_\varepsilon} \Vert_{C([0,T];\mathbb{L}^2(\Omega))} + \Vert \lambda_{h_\varepsilon}  - \lambda^{h_\varepsilon} \Vert_{C([0,T];\mathbb{L}^2(\Omega))} \leq C \Vert \phi_\varepsilon \Vert_{\mathbb{L}^2(\Omega)}.
	\end{equation}
\end{lemma}
\begin{proof}
\textbf{Step 1.} Putting $ \varphi_{h_\varepsilon} = u_{h_\varepsilon} - u^{h_\varepsilon} $, then $ \varphi \in C([0,T],H) $. It is enough to estimate $ \Vert \varphi_{h_\varepsilon} \Vert_{L^\infty(0,T;H)} $. From the unperturbed equations \eqref{VCHE_original} and the perturbed equations \eqref{VCHE_perturbed}, it follows that $ \varphi_{h_\varepsilon} $ is the solution of the equations
	\begin{equation} \label{estimate_varphi}
	\begin{cases}
		\partial_t(\varphi_{h_\varepsilon} +\alpha^2 A\varphi_{h_\varepsilon})+\nu A(\varphi_{h_\varepsilon}+\alpha^2 A\varphi_{h_\varepsilon})
		\\ \hspace{81pt}+\widetilde{B}(u_{h_\varepsilon},u_{h_\varepsilon}+\alpha^2 A u_{h_\varepsilon}) - \widetilde{B}(u^{h_\varepsilon},u^{h_\varepsilon}+\alpha^2 A u^{h_\varepsilon}) = 0,\\
		\varphi_{h_\varepsilon}(0)=\phi_0,
	\end{cases}
\end{equation}
where  \begin{align*} \widetilde{B} &(u_{h_\varepsilon},u_{h_\varepsilon}+\alpha^2 A u_{h_\varepsilon}) - \widetilde{B}(u^{h_\varepsilon},u^{h_\varepsilon}+\alpha^2 A u^{h_\varepsilon}) \\ &= \widetilde{B}(\varphi_{h_\varepsilon}, u_{h_\varepsilon}+\alpha^2 A u_{h_\varepsilon}) + \widetilde{B}(u^{h_\varepsilon}, \varphi_{h_\varepsilon}+ \alpha^2 A\varphi_{h_\varepsilon}) . \end{align*}
Taking the inner product of \eqref{estimate_varphi} with $ \varphi_{h_\varepsilon} $ yields the following equality
\begin{align*} \frac{1}{2} \frac{d}{dt}(|\varphi_{h_\varepsilon}|^2&+\alpha^2 \left \| \varphi_{h_\varepsilon} \right \|^2) +\nu (\left \| \varphi_{h_\varepsilon} \right \|^2 + \alpha^2 |A\varphi_{h_\varepsilon}|^2)\\& = -\widetilde{b}(u^{h_\varepsilon}, \varphi_{h_\varepsilon}+\alpha^2 A\varphi_{h_\varepsilon},\varphi_{h_\varepsilon}) .
\end{align*}
Using the estimate $ | \widetilde{b}(u^{h_\varepsilon}, \varphi_{h_\varepsilon}+\alpha^2 A\varphi_{h_\varepsilon},\varphi_{h_\varepsilon}) | \leq C \Vert \varphi_{h_\varepsilon} \Vert^2 + \nu \alpha^2 |A\varphi_{h_\varepsilon}|^2$, we obtain
\[ \frac{1}{2} \frac{d}{dt} \left(|\varphi_{h_\varepsilon}(t)|^2 + \alpha^2 \Vert \varphi_{h_\varepsilon} (t)\Vert^2\right) \leq C \Vert \varphi_{h_\varepsilon}(t) \Vert^2. \]
Thus, \[ |\varphi_{h_\varepsilon}(t)|^2 + \alpha^2 \Vert \varphi_{h_\varepsilon}(t) \Vert^2 \leq \Vert\varphi_{h_\varepsilon}(0)\Vert^2 + |\varphi_{h_\varepsilon}(0)|^2 + K.  \]
Now, by combining the above inequalites with the initial condition of $ \varphi_{h_\varepsilon} $, we get the desired estimate $ \Vert\varphi_{h_\varepsilon}\Vert_{C([0,T];\mathbb{L}^2(\Omega))} \leq \Vert \phi_\varepsilon\Vert_{\mathbb{L}^2(\Omega)}. $\\
\textbf{Step 2.} For $ \lambda_{h_\varepsilon}$ and $\lambda^{h_\varepsilon} $, we set $ \psi_{h_\varepsilon} =  \lambda_{h_\varepsilon} - \lambda^{h_\varepsilon} $. From the adjoint equation \eqref{adjoint eq}, we deduce that $ \psi_{h_\varepsilon} $ is the solution of the following equations
\begin{equation}\label{estimate_psi}
	\begin{cases}
		-\partial_t(\psi_{h_\varepsilon}+\alpha^2A\psi_{h_\varepsilon}) +\nu A(\psi_{h_\varepsilon}+\alpha^2A\psi_{h_\varepsilon})+\hat{B}(u_{h_\varepsilon},\lambda_{h_\varepsilon})-\hat{B}(u^{h_\varepsilon},\lambda^{h_\varepsilon}) = \\ \hspace{15pt} [L_u(x,t,u_{h_\varepsilon}(x,t),h_\varepsilon(x,t))-L_u(x,t,u^{h_\varepsilon},h_\varepsilon(x,t))]
		\hfill \text{ in } V' \text{ for a.e. } t\in [0,T],\\
		A\psi_{h_\varepsilon} = 0, \hspace{186.5pt} \text{ in } \partial \Omega \times [0,T],\\
		\psi_{h_\varepsilon}(T)+\alpha^2 A\psi_{h_\varepsilon}(T) = F_u(x,u_{h_\varepsilon}(x,T)) - F_u(x,u^{h_\varepsilon}(x,T)) \text{ in } V',
	\end{cases}
\end{equation}
Taking the inner product of \eqref{estimate_psi} with $ \psi_{h_\varepsilon} $ leads to
\begin{align}
	&-\frac{1}{2} \frac{d}{dt}  \left( |\psi_{h_\varepsilon}|^2 + \Vert \psi_{h_\varepsilon} \Vert^2 \right) + \nu \left( \Vert \psi_{h_\varepsilon} \Vert ^2 + \alpha^2 |A\psi_{h_\varepsilon}|^2 \right) + \nonumber \\
	& \widetilde{b}(\varphi_{h_\varepsilon}, \psi_{h_\varepsilon} + \alpha^2 A\psi_{h_\varepsilon},\lambda_{h_\varepsilon}) + \widetilde{b}(u^{h_\varepsilon}, \psi_{h_\varepsilon} +\alpha^2 A \psi_{h_\varepsilon},\psi_{h_\varepsilon}) +\widetilde{b}(\psi_{h_\varepsilon},\varphi_{h_\varepsilon} +\alpha^2 A \varphi_{h_\varepsilon},\lambda_{h_\varepsilon}) \nonumber\\&=L_u(x,t,u_{h_\varepsilon}(x,t),h_\varepsilon(x,t))-L_u(x,t,u^{h_\varepsilon},h_\varepsilon(x,t)). \label{innerproduct_psi}
\end{align}
Since $ u^{h_\varepsilon}, \lambda_{h_\varepsilon} \in L^{\infty}(0,T; \mathbb{H}_{\sigma}^2(\Omega) )$, utilizing Young's inequality together with the boundedness of $ \Vert \varphi_{h_\varepsilon} \Vert $ in Step 1, we obtain
\begin{align*}
	 |\widetilde{b}(u^{h_\varepsilon}, \psi_{h_\varepsilon} + \alpha^2 A \psi_{h_\varepsilon}, \psi_{h_\varepsilon})| & \leq C \Vert \psi_{h_\varepsilon} \Vert^2 + \frac{v \alpha^2}{2} |A \psi_{h_\varepsilon}|^2,\\
	|\widetilde{b}(\varphi_{h_\varepsilon}, \psi_{h_\varepsilon}+\alpha^2 A\psi_{h_\varepsilon}, \lambda_{h_\varepsilon})| &\leq |b(\varphi_{h_\varepsilon}, \psi_{h_\varepsilon} +\alpha^2 A \psi_{h_\varepsilon}, \lambda_{h_\varepsilon})| + |b(\lambda_{h_\varepsilon}, \psi_{h_\varepsilon}+ \alpha^2 A\psi_{h_\varepsilon}, \varphi_{h_\varepsilon})\\
	&\leq |b(\varphi_{h_\varepsilon}, \psi_{h_\varepsilon}, \lambda_{h_\varepsilon})| + |b(\varphi_{h_\varepsilon}, \alpha^2 A\varphi_{h_\varepsilon}, \lambda_{h_\varepsilon})|\\
	&\hspace{20pt} |b(\lambda_{h_\varepsilon}, \psi_{h_\varepsilon}, \varphi_{h_\varepsilon})| + |b(\lambda_{h_\varepsilon}, \alpha^2 A \psi_{h_\varepsilon}, \varphi_{h_\varepsilon})| \\
	& \leq C_1 \Vert \psi_{h_\varepsilon} \Vert^2 + \frac{\nu \alpha^2}{2} |A \psi_{h_\varepsilon}|^2,\\
	|\widetilde{b}(\psi_{h_\varepsilon}, \varphi_{h_\varepsilon} + \alpha^2 A \varphi_{h_\varepsilon}, \lambda_{h_\varepsilon})| &\leq C_2 \Vert \psi_{h_\varepsilon} \Vert^2.
\end{align*}
Substituting these inequalities in \eqref{innerproduct_psi} and integrating from $ t $ to $ T $ yield the following inequalities
\begin{align*}
	|\psi_{h_\varepsilon}(t)|^2 &+ \Vert \psi_{h_\varepsilon}(t) \Vert ^2 + C \int_{t}^{T} \Vert \psi_{h_\varepsilon}(s) \Vert \ ds \\
	& \leq |\psi_{h_\varepsilon}(T)|^2 + \alpha ^2 \Vert \psi_{h_\varepsilon}(T) \Vert^2 + K \\
	& \leq \left(F_u(x,u_{h_\varepsilon}(x,T)) - F_u(x,u^{h_\varepsilon}(x,T)), \psi_{h_\varepsilon}(T) \right),
\end{align*}
where we have used the Lipchitz condition \eqref{Liptchiz_Lu} of $ L_u $.  \\
Finally, the estimate $ \Vert \psi_{h_\varepsilon} \Vert_{C([0,T];\mathbb{L}^2(\Omega))} \leq C \Vert \phi_\varepsilon\Vert_{\mathbb{L}^2(\Omega)}  $ holds thanks to the Lipchitz assumption \eqref{Lipchitz_Fu} of $ F_u $ and the estimate of $ \Vert\varphi_{h_\varepsilon}\Vert_{C([0,T];\mathbb{L}^2(\Omega))} $ in Step 1.
\end{proof}

\begin{lemma} \label{transform_SSC}
	The second-order sufficient optimality condition \eqref{SSC} is equivalent to
	\begin{equation}
		\exists \, \tau > 0 \ \text{and} \ \exists \, \mu > 0 \ \text{such that} \ J''(\hat{h})\widetilde{h}^2 \geq \mu \Vert \widetilde{h} \Vert_{\mathbb{L}^2(Q)}^2 \ \forall \, \widetilde{h} \in E_{\hat{h}}^\tau,
	\end{equation}
where \begin{equation}
	E_{\hat{h}}^\tau = \{ \widetilde{h} \in \mathcal{T}_{\mathcal{H}}(\hat{h}) \ \text{and } \ J'(\hat{h})\widetilde{h} \leq \tau \Vert v \Vert_{\mathbb{L}^2(Q)} \}.
\end{equation}
\end{lemma}
\begin{proof}
	We argue by contradiction. Assume that there exists a sequence $ \left\{\left(h_k, \widetilde{h}_k\right) \right\}_{k = 0}^\infty$ in $ \mathbb{L}^2(Q) \times \mathbb{L}^2(Q) $ such that $ \displaystyle \|h_k - \hat{h}\|_{\mathbb{L}^2(Q)} \rightarrow 0, J''(h_k)\widetilde{h}_k \leq \frac{1}{k} \, \|\widetilde{h}_k\|_{\mathbb{L}^2(Q)} $, and $ \widetilde{h}_k \in E_{\hat{h}}^{1/k}  $.\\
	Putting $ \widetilde{h}_k = \dfrac{\widetilde{h}_k}{\|\widetilde{h}_k\|_{\mathbb{L}^2(Q)}} $, then we still have that $ \widetilde{h}_k \in E_{\hat{h}}^{1/k} $. Since $ \{\widetilde{h}_k\} $ is bounded in $ \mathcal{T}_{\mathcal{H}}(\hat{h}) $, selecting a subsequence if necessary, we obtain an element $ \widetilde{h} \in  \mathcal{T}_{\mathcal{H}}(\hat{h}) $ such that $ \widetilde{h}_k $ converges weakly to $ \widetilde{h} $ in $ \mathcal{T}_{\mathcal{H}}(\hat{h}) $. Furthermore, we get the inequalities $\displaystyle J'(\hat{h})\widetilde{h}_k \leq \frac{1}{k} $ and $\displaystyle J''(h_k)\widetilde{h}_k^2 \leq \frac{1}{k} $ for all $ k $.\\
	It is also clear from the properties of $ \mathcal{T}_{\mathcal{H}}(\hat{h}) $ that $ \widetilde{h} \in \mathcal{T}_{\mathcal{H}}(\hat{h}) $. Then arguing as Theorem \ref{SSC}, we infer the contradiction from $ \widetilde{h} \in C_0(\hat{h}) $ and $ J''(\hat{h})\widetilde{h} \leq 0 $.
\end{proof}

\begin{lemma}\label{lemma3}
	There exists $ \varepsilon_0 >0  $ such that
	\begin{equation}
		J''(\hat{h} + \theta(h_\varepsilon - \hat{h}))(h_\varepsilon - \hat{h})^2 \geq \frac{\mu}{2} \Vert h_\varepsilon - \hat{h} \Vert_{\mathbb{L}^2(Q)}^2 \ 
		\forall \, \varepsilon \in (0,\varepsilon_0) \ \text{and} \ \forall \, \theta \in (0,1),
	\end{equation}
where $ \mu $ is given by Lemma  \ref{transform_SSC}.
\end{lemma}

\begin{proof}
	Let us take $ \tau > 0 $ as in Lemma  \ref{transform_SSC}. We first prove that $ h_\varepsilon - \hat{h} $ belongs to $ E_{\hat{h}}^\tau $ for every sufficiently small $ \varepsilon $. Since $ h_\varepsilon - \hat{h} $ satisfies the sign condition defining $ \mathcal{T}_{\mathcal{H}}(\hat{h}) $, $ h_\varepsilon - \hat{h} $ belongs to $ \mathcal{T}_{\mathcal{H}}(\hat{h}) $. It remains to confirm that $ J'(\hat{h})(h_\varepsilon - \hat{h}) \leq \tau \|h_\varepsilon - \hat{h}\|_{\mathbb{L}^2(Q)} $ for $ \varepsilon $ small enough. \\
	To this end we set $\displaystyle v_\varepsilon = \frac{h_\varepsilon - \hat{h}}{\|h_\varepsilon - \hat{h}\|_{\mathbb{L}^2(Q)}} $. Then it holds $v_\varepsilon \in \mathcal{T}_{\mathcal{H}}(\hat{h}) $. Taking a subsequence, we can assume that $ v_\varepsilon \rightharpoonup v $ in $ \mathbb{L}^2(Q) $. In the proof of Theorem \ref{1st_theorem_stability}, it yields $ h_\varepsilon $ converges strongly to $ \hat{h} $ in $ \mathbb{L}^2(Q) $.\\
	Then we get
	\[ J'(\hat{h})v = \lim_{\varepsilon \rightarrow 0 }J'(\hat{h})v_\varepsilon \geq 0 \ \text{and} \ J'(\hat{h}) v = \lim_{\varepsilon \rightarrow 0} J'(h_\varepsilon)v_\varepsilon \leq 0. \]
	Hence, $\underset{\varepsilon \rightarrow 0}{\lim} J'(\hat{h})v_\varepsilon = J'(\hat{h})v = 0 $ holds. Since this happens for any convergent subsequence of $ \{v_\varepsilon\}_{\varepsilon>0} $, we infer that $ J'(\hat{h})v_\varepsilon \rightarrow 0 $ as $ \varepsilon \rightarrow 0 $ holds for the whole family, meaning that \[ \lim_{\varepsilon \rightarrow 0} J'(\hat{h})v_\varepsilon = 0. \]
	Therefore, there exists $ \varepsilon_0 $ such that \[ J'(\hat{h})v_\varepsilon \leq \tau \ \text{for every} \ \varepsilon \in (0,\varepsilon_0),\]
	or equivalently,
	\[ J'(\hat{h})(h_\varepsilon - \hat{h}) \leq \tau \|h_\varepsilon - \hat{h}\|_{\mathbb{L}^2(Q)}, \]
	which implies that $ h_\varepsilon - \hat{h} $ belongs to $ E_{\hat{h}}^\tau $.\\
	Then, since $ \hat{h} +\theta(h_\varepsilon - \hat{h}) $ converges strongly to $ \hat{h} $ in $ \mathbb{L}^2(Q) $ and $ J $ is of class $ C^2 $, we estimate $ \underset{\varepsilon \rightarrow 0}\lim J''(\hat{h} + \theta (h_\varepsilon - \hat{h})) (h_\varepsilon - \hat{h})^2 = J''(\hat{h})(h_\varepsilon - \hat{h})^2 $. Thus, we infer $\displaystyle J''(\hat{h} + \theta (h_\varepsilon - \hat{h}))(h_\varepsilon - \hat{h})^2 \geq \frac{\mu}{2}\|h_\varepsilon - \hat{h}\|_{\mathbb{L}^2(Q)}^2. $
\end{proof}

We are now ready to prove Theorem \ref{3rd_theorem_stability}.

\begin{proof}[Proof of Theorem \ref{3rd_theorem_stability}] 
	The local optimality of $ \hat{h} $ and $ h_\varepsilon $ leads to 
	\[ J'(\hat{h})(h_\varepsilon - \hat{h}) \geq 0 \ \text{and} \ J'_\varepsilon(h_\varepsilon)(\hat{h} - h_\varepsilon) \geq 0. \]
	Hence, \[ [J'(h_\varepsilon) - J'(\hat{h})](h_\varepsilon - \hat{h}) \leq [J'_\varepsilon(h_\varepsilon) - J'(h_\varepsilon)](\hat{h} - h_\varepsilon). \]
	Then, invoking from the mean value theorem, Lemma \ref{lemma1} and Lemma \ref{lemma3}, we deduce for $ \varepsilon $ small enough,
	\begin{align*}
		\frac{\mu}{2}&\|h_\varepsilon-\hat{h}\|_{\mathbb{L}^2(Q)}^2 \leq J''(\hat{h} +\theta_k(h_\varepsilon - \hat{h}))(h_\varepsilon-\hat{h})^2\\
		&=[J'(h_\varepsilon) - J'(\hat{h})](h_\varepsilon-\hat{h}) \leq [J'_\varepsilon(h_\varepsilon) - J'(h_\varepsilon)](\hat{h} - h_\varepsilon)\\
		&=(\lambda^{h_\varepsilon} - \lambda_{h_\varepsilon}, \hat{h} - h_\varepsilon)_Q \leq \|\lambda^{h_\varepsilon} - \lambda_{h_\varepsilon}\|_{\mathbb{L}^2(Q)} \|\hat{h} - h_\varepsilon\|_{\mathbb{L}^2(Q)}\\
		&\leq C\Vert \phi_\varepsilon \Vert_{\mathbb{L}^2(\Omega)} \|\hat{h} - h_\varepsilon\|_{\mathbb{L}^2(Q)}.
	\end{align*}
	This completes the proof.
\end{proof}

	\vskip 0.5cm
	
	\noindent{\bf Acknowledgements.} This research is funded by Vietnam National Foundation for Science and Technology Development (NAFOSTED) under grant number 101.02--2023.29. The authors also would like to thank Hanoi National University of Education for providing a fruitful working environment.
	
	\vskip 0.5cm 
	
\noindent {\bf Disclosure statement.}  The authors declare that they have no potential conflict of interest.

\end{document}